\documentclass[10pt,reqno]{amsart}

\usepackage{enumerate}
\usepackage{amsmath, amssymb, amsthm,amsfonts}
\usepackage{mathrsfs}
\usepackage{esint}
\usepackage{xcolor}
\usepackage{mathtools}
\usepackage{hyperref}
\usepackage{bm}
\usepackage{accents}
\usepackage[foot]{amsaddr}
\usepackage{todonotes}

\makeatletter

\numberwithin{equation}{section}

\newtheorem{thm}{Theorem}[section]
\newtheorem{theorem}[thm]{Theorem}
\newtheorem{lemma}[thm]{Lemma}
\newtheorem{corollary}[thm]{Corollary}

\theoremstyle{definition}

\newtheorem{assumption}[thm]{Assumption}
\theoremstyle{remark}

\newtheorem{remark}[thm]{\bf{Remark}}

\newcommand\aint{-\hspace{-0.38cm}\int}

\newcommand\bM{\mathbb{M}}
\newcommand\bN{\mathbb{N}}

\newcommand\bR{\mathbb{R}}
\newcommand\bS{\mathbb{S}}

\newcommand\bZ{\mathbb{Z}}

\newcommand\cF{\mathcal{F}}

\newcommand\cS{\mathcal{S}}

\begin{document}

\title[]{Sobolev regularity for nondivergence form nonlocal equations with irregular kernels}

\author{Junhee Ryu}
\address{School of Mathematics, Korea Institute for Advanced Study, 85 Hoegi-ro, Dongdaemun-gu, Seoul, 02455, Republic of Korea}
\email{junhryu@kias.re.kr}
\thanks{J. Ryu was supported by a KIAS Individual Grant (MG101501) at Korea Institute for Advanced Study and by the National Research Foundation of Korea (NRF) grant funded by the Korea government (MSIT) (No. RS-2026-25477288).}



\begin{abstract}
We consider the nonlocal elliptic equation $L_Ku-\lambda u=f$ in $\mathbb{R}^d$, where $L_K$ is a nondivergence form operator of order $\sigma\in(0,2)$ with kernel $K(x,y)$.
We assume that $K(x,y)$ is Dini continuous in $x$ uniformly in $y$ and merely measurable in $y$.
Under the usual upper bound assumption, we prove the boundedness of $L_K$ from $H_p^\sigma(\mathbb{R}^d)$ to $L_p(\mathbb{R}^d)$ for any $p\in(1,\infty)$.
We also establish maximal $L_p$-regularity and unique solvability for all sufficiently large $\lambda$, provided that an additional nondegeneracy condition holds.
We further extend these results to weighted spaces with arbitrary Muckenhoupt $A_p$ weights and obtain estimates that are robust as $\sigma\to2$.
We also construct a kernel that is jointly uniformly continuous in $(x,y)$, for which both operator boundedness and the maximal regularity estimate fail.
This is in sharp contrast to the classical $W_p^2$ theory for second-order elliptic equations in nondivergence form, which allows VMO coefficients.
\end{abstract}

\maketitle

\section{Introduction}

In this paper, we consider the nonlocal elliptic equation
\begin{equation} \label{intro_equation}
L_Ku-\lambda u=f \quad \text{in } \bR^d,
\end{equation}
where $d\geq1$, $\lambda\geq0$, and $L_K$ is a nondivergence form operator of order $\sigma\in(0,2)$ defined by
\begin{equation} \label{intro_op}
L_Ku(x):=\int_{\bR^d}\delta_yu(x)K(x,y)\frac{dy}{|y|^{d+\sigma}}.
\end{equation}
Here,
\begin{equation} \label{eq_differ}
y^{(\sigma)} := y1_{\{\sigma>1\}} + y1_{\{\sigma=1\}}1_{\{|y|<1\}}, \quad \delta_yu(x) := u(x+y)-u(x)-y^{(\sigma)}\cdot\nabla u(x).
\end{equation}
The kernel $K(x,y)$ is allowed to be merely measurable in the jump variable $y$.

From a probabilistic viewpoint, $L_K$ has the form of the infinitesimal generator of a pure jump L\'evy-type process, with $K(x,y)|y|^{-d-\sigma}dy$ describing the intensity of jumps from $x$ to $x+y$.
More precisely, for a Borel set $B$ bounded away from the origin, the integral $\int_B K(x,y)|y|^{-d-\sigma}dy$ represents the instantaneous rate of jumps from $x$ to $x+B$.
In the case of L\'evy processes, the jump intensities are independent of the starting points and depend only on the jump displacements.
The spatial dependence of $K(x,y)$ allows L\'evy-type processes to describe jumps whose intensities vary with their starting points.
When $K$ is independent of $x$, one recovers a class of L\'evy processes.

We establish maximal $L_p$-regularity and solvability for \eqref{intro_equation} under Dini continuity of the kernel in the spatial variable.
More precisely, we assume
\begin{equation} \label{intro_dini}
|K(x_1,y)-K(x_2,y)| \leq (2-\sigma) \omega(|x_1-x_2|), \quad 0<|x_1-x_2|\leq1,
\end{equation}
where $\omega$ is a nondecreasing function such that
\begin{equation*}
N_\omega := \int_0^1\omega(r)\frac{dr}{r}<\infty.
\end{equation*}
We emphasize that no regularity assumption on $K(x,y)$ in the $y$-variable is imposed, and only measurability is assumed.
Under the upper bound and nondegeneracy conditions (see Assumptions \ref{as_bound} and \ref{as_nonde}), we prove that, for every $p\in(1,\infty)$ and $f\in L_p(\bR^d)$, equation \eqref{intro_equation} has a unique solution $u\in H_p^\sigma(\bR^d)$ for all sufficiently large $\lambda$, and
\begin{equation} \label{intro_est}
\|(-\Delta)^{\sigma/2}u\|_{L_p(\bR^d)}+\lambda\|u\|_{L_p(\bR^d)}\leq N\|f\|_{L_p(\bR^d)}.
\end{equation}
The Dini condition \eqref{intro_dini} cannot, in general, be relaxed to mere uniform continuity.
Indeed, in \cite{L26}, it was shown that, for every $\sigma\in(0,2)$ and $p\in[1,\infty)$, the interior $L_p$-estimate for nonlocal equations can fail even when the kernel is jointly uniformly continuous in $(x,y)$.
Adapting this counterexample, we obtain a kernel that is jointly uniformly continuous in $(x,y)$ and for which the global estimate \eqref{intro_est} also fails.

A key difference between nonlocal equations and second-order elliptic equations lies in the spatial regularity assumptions on their coefficients.
For second-order uniformly elliptic equations in nondivergence form, global $W_p^2$ estimates and solvability are available under VMO assumptions on the leading coefficients.
In \cite{CFL93}, the authors established global estimates and solvability for the Dirichlet problem, while \cite{Kry07} developed a whole space $L_p$ theory for elliptic and parabolic equations.
For further results on broader classes of coefficients, including partially small mean oscillation conditions, see \cite{KK07,K09,DK18}.
These results illustrate that Sobolev estimates for second-order equations can be obtained even for a substantial class of discontinuous coefficients.
We also refer the reader to \cite{AQ02,BW04,DK10,DK11} for results on second-order equations in divergence form.

For nondivergence form nonlocal equations, $L_p$ estimates and solvability were established in \cite{DK12} when the kernel is independent of $x$ and merely measurable in $y$.
The theory has also been developed for parabolic equations with general L\'evy measures, including singular measures, in \cite{Z13,DR25}.
Unlike measures of the form $K(y)|y|^{-d-\sigma}dy$, the L\'evy measures $\nu(dy)$ considered in these works are not necessarily absolutely continuous with respect to the Lebesgue measure.

For spatially dependent kernels, parabolic Sobolev solvability under H\"older type assumptions in $x$ was obtained in \cite{MP14} for sufficiently large $p$.
In \cite{DJK23}, the operator boundedness was established for all $p\in(1,\infty)$ under a H\"older continuity assumption in $x$.
We also note that unweighted elliptic and parabolic solvability was proved in the same paper.
Further results for parabolic equations were obtained in \cite{DL23}, including Sobolev estimates and solvability for equations with nonlocal time derivatives, as well as weighted mixed-norm estimates and solvability for equations with the ordinary time derivative.
Parabolic solvability was obtained in \cite{Z13_var} under Dini-type regularity assumptions in both $x$ and $y$.
The additional Dini condition in $y$ near the origin plays a crucial role in the perturbation argument.

The main contribution of this paper is to extend the weighted operator estimates of \cite{DJK23} to the weaker Dini continuity condition \eqref{intro_dini}, without imposing any additional regularity assumptions in the jump variable.
We also establish maximal $L_p$-regularity and unique solvability for \eqref{intro_equation} in weighted Bessel potential spaces for arbitrary Muckenhoupt $A_p$ weights.

We next discuss regularity results for nonlocal equations in divergence form, where weaker assumptions on the kernels can be imposed.
To describe this structure, let $A:\bR^d\times\bR^d\to\bR$ be measurable and satisfy $A(x,z)=A(z,x)$ and $0<\nu\leq A(x,z)\leq\Lambda$.
Formally, we write
\begin{equation*}
L_Au(x) := \text{p.v.} \int_{\bR^d}(u(x)-u(z))A(x,z) \frac{dz}{|x-z|^{d+\sigma}}.
\end{equation*}
For $f\in L_{1,loc}(\bR^d)$, we say that $u\in H_2^{\sigma/2}(\bR^d)$ is a weak solution of $L_Au=f$ in $\bR^d$ if
\begin{equation*}
\frac{1}{2} \int_{\bR^d\times\bR^d} \frac{(u(x)-u(z))(\varphi(x)-\varphi(z))A(x,z)}{|x-z|^{d+\sigma}} dzdx = \int_{\bR^d} f(x)\varphi(x) dx
\end{equation*}
for any $\varphi\in C_c^\infty(\bR^d)$.

In this weak framework, higher Sobolev regularity is available when the coefficient $A(x,z)$ is VMO in $(x,z)$. See \cite{N23VMO,N23high}.
More precisely, under this assumption, it was shown in \cite{N23high} that if $d>\sigma$ and $2<p<\infty$, then a weak solution $u\in H_2^{\sigma/2}(\bR^d)$ with $f\in L_{p,loc}(\bR^d)$ satisfies
\begin{equation*}
u\in W_{loc}^{t,p}(\bR^d), \quad \frac{\sigma}{2}\leq t<\min\{\sigma,1\}.
\end{equation*}
For earlier results on $L_p$ regularity theory, we refer the reader to \cite{KMS15,ABES19,MSY21}.
For related developments under VMO or SMO (small BMO) assumptions, including fractional $p$-Laplacian type equations, see \cite{BK25,DN25}.
We also refer the reader to \cite{Bass09,CS09,DK13,SS16,DZ18,DyK20,KW24} for related H\"older and Schauder regularity theories.

Although these results establish higher Sobolev regularity under VMO assumptions, they concern divergence form nonlocal equations and do not directly extend to the nondivergence form setting.
In particular, they do not address maximal $L_p$-regularity for the operators considered in this paper.

The boundedness of nonlocal operators presents an additional difficulty.
For a second-order operator $a^{ij}(x)D_{ij}$, it is straightforward to obtain boundedness from $W_p^2$ to $L_p$ when the coefficients are bounded.
In contrast, for the nonlocal operator \eqref{intro_op}, the boundedness from $H_p^\sigma$ to $L_p$ does not follow directly from the boundedness of the kernel due to its integral structure.
To address this difficulty, we use the boundedness estimates for spatially independent kernels together with the Dini continuity of the kernel in the spatial variable.

We now briefly describe the idea of the proof.
We first freeze the spatial variable of the kernel on a grid of mesh size $r$, that is, we replace $K(x,y)$ by $K(x_m,y)$ when $x$ is near a grid point $x_m$.
However, unlike the case of local operators, the error of this approximation cannot be controlled simply by the smallness of $K(x,y)-K(x_m,y)$, since the operator involves jumps of all sizes.
To address this difficulty, we also decompose the integral with respect to $y$ according to the size of the jumps, and choose the mesh size of the grid depending on the size of the jumps.
The Dini condition allows us to sum the resulting estimates over all sizes of jumps.
As a result, in sufficiently small balls, $L_K$ is a small perturbation of a spatially independent operator, and the main results follow from a standard localization argument and the method of continuity.

The paper is organized as follows.
In Section \ref{sec_main}, we introduce the function spaces, assumptions, and main results.
In Section \ref{sec_aux}, we collect auxiliary results, and in Section \ref{sec_analysis}, we develop the estimates for nonlocal operators.
In Section \ref{sec_proof}, we prove operator boundedness, maximal regularity, and solvability under the Dini condition.
In Section \ref{sec_counter}, we construct global counterexamples with kernels that are jointly uniformly continuous in $(x,y)$.

\section{Main results} \label{sec_main}

\subsection{Notation and function spaces}

In this subsection, we introduce notation and function spaces used throughout the paper.
 We use ``$:=$'' or ``$=:$'' to denote a definition.
  We denote the set of natural numbers by $\bN$, and write $\bN_0:=\bN\cup\{0\}$. As usual, $\bR^d$ stands for the Euclidean space of points $x=(x_1,\dots,x_d)$, and we write $\bR:=\bR^1$.
  For $r>0$ and $x\in \bR^d$, we denote by $B_r(x) := \{y\in \bR^d: |x-y|<r\}$ the open ball centered at $x$ with radius $r$. When $x=0$, we write $B_r:=B_r(0)$.
We use $D^n_x u$ to denote the partial derivatives of order $n\in\bN_0$ with respect to the spatial variables, and $D_xu:=\nabla u:= D_x^1u$. 
For $E\subset \bR^{d}$ such that $0<|E|<\infty$, we denote by
\begin{equation*}
\aint_{E} f(x) dx := \frac{1}{|E|}\int_E f(x) dx
\end{equation*}
the average of $f$ over a measurable set $E$.
For $D\subset \bR^{d}$, we denote by $C_c^\infty(D)$ the collection of all infinitely differentiable functions with compact support in $D$.
We denote the Schwartz space and the space of tempered distributions by $\cS(\bR^d)$ and $\cS'(\bR^d)$, respectively. 
We define the Fourier transform and its inverse by
\begin{equation*}
\begin{gathered}
\hat{f}(\xi) := \cF[f](\xi) := \int_{\bR^d} e^{-ix\cdot\xi} f(x) dx,
\\
\check{f}(x) := \cF^{-1}[f](x) := (2\pi)^{-d} \int_{\bR^d} e^{ix\cdot\xi} f(\xi) d\xi.
\end{gathered}
\end{equation*}

We now introduce the function spaces and related notation.
For $p\in(1,\infty)$ and a locally integrable nonnegative weight $w(x)$ on $\bR^d$, we define $L_{p,w} := L_p(\bR^d, w dx)$, equipped with the norm
\begin{equation*}
    \left\|u\right\|_{L_{p,w}} = \left( \int_{\bR^d} |u(x)|^p w(x) dx \right)^{1/p}.
\end{equation*}
We say that $w$ is in the $A_p$ Muckenhoupt class of weights if
\begin{equation*}
    [w]_{A_p} := \sup_{r>0,x_0\in\bR^d} \left(\aint_{B_r(x_0)} w(x) dx\right) \left(\aint_{B_r(x_0)} (w(x))^{1/(1-p)} dx\right)^{p-1} < \infty.
\end{equation*}
The condition $w\in A_p$ is equivalent to the boundedness of the Hardy-Littlewood maximal operator $\bM$ on $L_{p,w}$, where
\begin{equation*}
\bM u(x) := \sup_{r>0} \aint_{B_r(x)}|u(y)| dy.
\end{equation*}
Indeed, for $w\in A_p$ with $[w]_{A_p}\leq K_0$, by the weighted maximal function theorem, we have
\begin{equation} \label{eqmaximal}
\|\bM u\|_{L_{p,w}} \leq N\|u\|_{L_{p,w}},
\end{equation}
where $N$ depends only on $d,p,K_0$ (see e.g. \cite[Theorem 2.2]{DK18}).

For $\sigma\in(0,2)$, $p\in(1,\infty)$, and $w\in A_p$, we denote by $H_{p,w}^\sigma$ the class of tempered distributions satisfying
\begin{equation*}
\left\|u\right\|_{H_{p,w}^\sigma} = \left\|(1-\Delta)^{\sigma/2}u\right\|_{L_{p,w}} < \infty.
\end{equation*}
Here, the operator $(1-\Delta)^{\sigma/2}$ is defined by
\begin{equation*}
    (1-\Delta)^{\sigma/2}u := \cF^{-1}\left\{(1+|\cdot|^2)^{\sigma/2}\cF[u]\right\}.
\end{equation*}
In the particular case when $w=1$, we write $H_p^{\sigma} :=H_{p,1}^{\sigma}$ and $L_p:=L_{p,1}$.
Moreover, for $[w]_{A_p}\leq K_0$,
\begin{equation} \label{eq10310040}
\|u\|_{H_{p,w}^\sigma} \approx \left(\|u\|_{L_{p,w}} + \|(-\Delta)^{\sigma/2}u\|_{L_{p,w}}\right),
\end{equation}
where $(-\Delta)^{\sigma/2}$ denotes the fractional Laplacian.
Here, the constants in the norm equivalence depend only on $d,p$, and $K_0$. See \cite[Lemma 3.4]{DJK23}.

 \begin{remark} \label{rem_dense}
The space $C_c^\infty(\bR^d)$ is dense in $H^\sigma_{p,w}$.
     Indeed, let $u\in H_{p,w}^\sigma$ and take a nonnegative function $\varphi\in C_c^\infty(B_1)$ such that $\int_{\bR^d}\varphi(x) dx=1$. We consider the mollification
\begin{equation*}
u^\varepsilon := u*\varphi_\varepsilon,
\end{equation*}
where $\varphi_\varepsilon(x) := \varepsilon^{-d}\varphi(x/\varepsilon)$.
Then we see that $u^\varepsilon\in C^\infty(\bR^d)$ and
\begin{equation} \label{eq10271414}
(-\Delta)^{\sigma/2}u^\varepsilon = \left((-\Delta)^{\sigma/2}u\right)*\varphi_\varepsilon.
\end{equation}
By the definition of convolution,
\begin{equation*}
\left|u^\varepsilon(x)\right| \leq N \aint_{B_\varepsilon(x)}|u(y)|dy \leq N\bM u(x).
\end{equation*}
Since $u*\varphi_\varepsilon\to u$ almost everywhere as $\varepsilon\to0$, by \eqref{eqmaximal} and the dominated convergence theorem, $u*\varphi_\varepsilon\to u$ in $L_{p,w}$.
Due to \eqref{eq10271414}, one can also derive that $(-\Delta)^{\sigma/2}u^\varepsilon\to (-\Delta)^{\sigma/2}u$ in $L_{p,w}$ as $\varepsilon\to0$.
Consequently, $u^\varepsilon\to u$ in $H_{p,w}^\sigma$ as $\varepsilon\to0$.

Now we show that $u^\varepsilon$ can be approximated by smooth compactly supported functions. For each fixed $\varepsilon>0$, its derivatives up to order two belong to $L_{p,w}$. By \cite[Theorem 1.2]{SV20},
\begin{equation} \label{eq11161331}
\|v\|_{H_{p,w}^\sigma} \leq N \left(\|v\|_{L_{p,w}} + \|\Delta v\|_{L_{p,w}}\right).
\end{equation}
Let $\eta_n\in C_c^\infty(B_{2n})$ be cutoff functions such that $0\leq\eta_n\leq1$, $\eta_n=1$ on $B_n$, and $|D_x^k\eta_n|\leq N n^{-k}$ for $k=1,2$. Then one can see that
\begin{equation*}
\left\|u^\varepsilon\eta_n-u^\varepsilon\right\|_{L_{p,w}} + \left\|\Delta(u^\varepsilon\eta_n-u^\varepsilon)\right\|_{L_{p,w}}\to0 \quad \text{ as } n\to\infty,
\end{equation*}
which together with \eqref{eq11161331} yields $u^\varepsilon\eta_n\to u^\varepsilon$ in $H_{p,w}^\sigma$. Since $u^\varepsilon\eta_n\in C_c^\infty(\bR^d)$, we obtain the desired result.
 \end{remark}

\subsection{Nonlocal equations with Dini continuous kernels}

Throughout this paper, we assume that the kernel $K(x,y)$ in \eqref{intro_op} is measurable in $(x,y)$.
We now state the assumptions on the kernel $K(x,y)$. First, we introduce the upper bound and cancellation condition.

\begin{assumption} \label{as_bound}
There is a constant $\Lambda>0$ such that
\begin{equation} \label{eq_upper}
0\leq K(x,y)\leq (2-\sigma)\Lambda.
\end{equation}
If $\sigma=1$, then, for any $x\in\bR^d$,
\begin{equation} \label{eq_cancel}
\int_{r<|y|<R} yK(x,y)\frac{dy}{|y|^{d+1}}=0, \quad 0<r<R<\infty.
\end{equation}
\end{assumption}

To establish solvability of the equation, we impose the following nondegeneracy condition.

\begin{assumption} \label{as_nonde}
There is a constant $\nu>0$ such that, for every $x,\xi\in\bR^d$,
\begin{equation} \label{eq_nonde}
\int_{|\xi\cdot y|\leq1} |\xi\cdot y|^2K(x,y) \frac{dy}{|y|^{d+\sigma}} \geq \nu |\xi|^\sigma.
\end{equation}
\end{assumption}

We now state our Dini continuity assumption on $K$.

\begin{assumption} \label{as_dini}
For any $0<|x_1-x_2|\leq1$ and $y\in\bR^d$,
\begin{equation} \label{eq_dini}
|K(x_1,y)-K(x_2,y)| \leq (2-\sigma) \omega(|x_1-x_2|),
\end{equation}
where $\omega:(0,1]\to[0,\infty)$ is a nondecreasing function such that
\begin{equation*}
    \int_0^1\omega(r)\frac{dr}{r}<\infty.
\end{equation*}
\end{assumption}

To describe the dependence of the constants in our main results, we introduce some notation.
Let 
\begin{equation*}
N_\omega := \int_0^1\omega(r)\frac{dr}{r}<\infty.
\end{equation*}
To specify the dependence of the constants uniformly as $\sigma\to2$, we also define
\begin{equation} \label{eq_sigma}
\widetilde{\sigma} :=
\begin{cases}
(\sigma_0,\sigma_1), &\text{ when } 0<\sigma_0\leq\sigma\leq\sigma_1<1,
\\
1, &\text{ when } \sigma=1,
\\
\sigma_0, &\text{ when } 1<\sigma_0\leq\sigma<2.
\end{cases}
\end{equation}

Our first result is about the boundedness of the operator.

\begin{theorem} \label{thm_bound}
Let $d\geq1$, $\sigma\in(0,2)$, $p\in(1,\infty)$, $K_0>0$, and $w\in A_p$ with $[w]_{A_p}\leq K_0$. Suppose that $K$ satisfies Assumptions \ref{as_bound} and \ref{as_dini}. Then $L_K$ is a continuous operator from $H_{p,w}^\sigma$ to $L_{p,w}$. Moreover,
\begin{equation} \label{est_bound}
\|L_Ku\|_{L_{p,w}} \leq N (\Lambda+N_\omega) \|u\|_{H_{p,w}^\sigma},
\end{equation}
where $N=N(d,p,K_0,\widetilde{\sigma})$ and $\widetilde{\sigma}$ is given by \eqref{eq_sigma}.
\end{theorem}

Next, we introduce the solvability result.

\begin{theorem} \label{thm_dini}
Let $d\geq1$, $\sigma\in(0,2)$, $p\in(1,\infty)$, $K_0>0$, and $w\in A_p$ with $[w]_{A_p}\leq K_0$.
Suppose that $K$ satisfies Assumptions \ref{as_bound}, \ref{as_nonde}, and \ref{as_dini}. Then there is $\lambda_0=\lambda_0(d,p,K_0,\nu,\Lambda,\widetilde{\sigma},\omega)>0$ such that for any $\lambda\geq\lambda_0$ and $u\in H_{p,w}^\sigma$ satisfying
\begin{equation} \label{eq_thm}
L_Ku-\lambda u=f,
\end{equation}
we have
\begin{equation} \label{est_main}
\|(-\Delta)^{\sigma/2}u\|_{L_{p,w}} + \lambda \|u\|_{L_{p,w}} \leq N \|f\|_{L_{p,w}},
\end{equation}
where $N=N(d,p,K_0,\nu,\Lambda,\widetilde{\sigma})$ and $\widetilde{\sigma}$ is defined in \eqref{eq_sigma}.
Moreover, for any $f\in L_{p,w}$ and $\lambda\geq\lambda_0$, there exists a unique solution $u\in H_{p,w}^\sigma$ to \eqref{eq_thm}. 
\end{theorem}

\begin{remark}
    $(i)$ The boundedness of the operator holds without imposing the nondegeneracy condition \eqref{eq_nonde}.

    $(ii)$ In both Theorems \ref{thm_bound} and \ref{thm_dini}, the estimates \eqref{est_bound} and \eqref{est_main} are robust in the sense that the corresponding constants $N$ do not blow up as $\sigma\to2$.
\end{remark}

\subsection{Failure of \texorpdfstring{$H_p^\sigma$}{Hpsigma} regularity under uniform continuity}

We next demonstrate the failure of $H_p^\sigma$ regularity
under uniform continuity, even when the continuity condition
is imposed jointly in $(x,y)$.

We denote by $\text{BUC}(\bR^d\times\bR^d)$ the space of bounded functions that are uniformly continuous jointly in $(x,y)$.

\begin{theorem} \label{thm_counter}
Let $d\geq1$, $p\in(1,\infty)$, and $\sigma\in(0,2)$. There exists a kernel $K\in\text{BUC}(\bR^d\times\bR^d)$ satisfying Assumptions \ref{as_bound} and \ref{as_nonde} such that the following assertions hold.

$(i)$ There is no finite constant $N$ such that for every $u\in C_c^\infty(\bR^d)$,
\begin{equation} \label{eq_buc_bound}
\|L_Ku\|_{L_p}\leq N\|u\|_{H_p^\sigma}.
\end{equation}

$(ii)$ For every fixed $\lambda\geq0$, there is no finite constant $N$ such that for every $u\in C_c^\infty(\bR^d)$,
\begin{equation} \label{eq_buc_est}
\|(-\Delta)^{\sigma/2}u\|_{L_p} \leq N \left(\|L_Ku-\lambda u\|_{L_p}+\|u\|_{L_p}\right).
\end{equation}
\end{theorem}

\section{Auxiliary results} \label{sec_aux}

In this section, we present some auxiliary results.

\begin{lemma} \label{lem11041430}
For any $\sigma\in(0,1)$ and $u\in\cS(\bR^d)$,
\begin{equation*}
u(x) = N_0 \int_{\bR^d} |z|^{-d+\sigma} (-\Delta)^{\sigma/2} u(x-z) dz
\end{equation*}
and
\begin{equation*}
u(x+y)-u(x) = N_0 \int_{\bR^d} \left( |z+y|^{-d+\sigma}-|z|^{-d+\sigma} \right) (-\Delta)^{\sigma/2} u(x-z) dz,
\end{equation*}
where $N_0 := \frac{\Gamma((d-\sigma)/2)}{2^\sigma\pi^{d/2}\Gamma(\sigma/2)}$.
\end{lemma}

\begin{proof}
See \cite[Lemma 3.1]{DJK23} or \cite[Lemma 2.2]{K84}.
\end{proof}

In the following lemma, we write $u\in L_{p,loc}$ if $\eta u\in L_p$ for any $\eta\in C_c^\infty(\bR^d)$.

\begin{lemma}
Let $\gamma\in(0,2)$, $R\in(0,\infty)$, and $u\in L_{1,loc}$. Then we have
\begin{equation} \label{eq10291650}
\int_{B_R^c} |u(x+y)| \frac{dy}{|y|^{d+\gamma}} \leq N R^{-\gamma} \gamma^{-1} \bM u(x)
\end{equation}
and
\begin{equation} \label{eq10291653}
\int_{B_R} |u(x+y)| \frac{dy}{|y|^{d-\gamma}} \leq N R^\gamma \gamma^{-1} \bM u(x),
\end{equation}
where $N$ depends only on $d$.
\end{lemma}

\begin{proof}
 Using $B_R^c=\cup_{j=0}^{\infty} B_{2^{j+1}R}\setminus B_{2^jR}$, we see that
\begin{align*}
\int_{B_R^c} |u(x+y)| \frac{dy}{|y|^{d+\gamma}} &= \sum_{j=0}^\infty \int_{B_{2^{j+1}R}\setminus B_{2^jR}} |u(x+y)| \frac{dy}{|y|^{d+\gamma}}
\\
&\leq \sum_{j=0}^\infty R^{-(d+\gamma)} 2^{-j(d+\gamma)} \int_{B_{2^{j+1}R}\setminus B_{2^jR}}|u(x+y)| dy
\\
&\leq NR^{-\gamma} \sum_{j=0}^\infty 2^{-\gamma j} \aint_{B_{2^{j+1}R}} |u(x+y)| dy
\\
&\leq NR^{-\gamma} (1-2^{-\gamma})^{-1} \bM u(x) \leq NR^{-\gamma}\gamma^{-1} \bM u(x),
\end{align*}
which proves \eqref{eq10291650}. Similarly, for \eqref{eq10291653},
\begin{align*}
\int_{B_R} |u(x+y)| \frac{dy}{|y|^{d-\gamma}} &= \sum_{j=0}^\infty \int_{B_{2^{-j}R}\setminus B_{2^{-j-1}R}} |u(x+y)| \frac{dy}{|y|^{d-\gamma}}
\\
&\leq 2^d \sum_{j=0}^\infty R^{-(d-\gamma)} 2^{j(d-\gamma)} \int_{B_{2^{-j}R}\setminus B_{2^{-j-1}R}} |u(x+y)| dy
\\
&\leq NR^\gamma \sum_{j=0}^\infty 2^{-j\gamma} \bM u(x)
\\
&= NR^\gamma (1-2^{-\gamma})^{-1} \bM u(x) \leq NR^\gamma \gamma^{-1} \bM u(x).
\end{align*}
Here, in the first inequality, we used $|y|^{\gamma-d}\leq (2^{-j-1}R)^{\gamma-d}$ for $\gamma\leq d$ and $|y|^{\gamma-d}\leq (2^{-j}R)^{\gamma-d}$ for $\gamma>d$.
The lemma is proved.
\end{proof}

\begin{lemma}
Let $c>0$ and $j_0\in\bZ$ be such that $c2^{1-j_0} \leq 1$. Suppose that $\omega:(0,1]\to[0,\infty)$ is a nondecreasing function. Then
\begin{equation} \label{eq10301427}
\sum_{j\geq j_0}\omega(c2^{-j}) \leq \frac{1}{\log2} N_\omega(c2^{1-j_0}),
\end{equation}
where
\begin{equation*}
    N_\omega(R) := \int_0^R\omega(r)\frac{dr}{r}.
\end{equation*}
\end{lemma}

\begin{proof}
Since $\omega$ is nondecreasing,
\begin{equation*}
\int_{c2^{-j}}^{c2^{1-j}} \omega(r) \frac{dr}{r} \geq \omega(c2^{-j}) \int_{c2^{-j}}^{c2^{1-j}} \frac{dr}{r} = \omega(c2^{-j})\log{2}.
\end{equation*}
Thus,
\begin{equation*}
\sum_{j\geq j_0} \omega(c2^{-j}) \leq \frac{1}{\log{2}} \int_{0}^{c2^{1-j_0}} \omega(r)\frac{dr}{r} = \frac{1}{\log2} N_\omega(c2^{1-j_0}),
\end{equation*}
which proves \eqref{eq10301427}.
The lemma is proved.
\end{proof}

\begin{lemma}
Let $\sigma\in(0,2)$, $\sigma'\in(0,\sigma)$, $p\in(1,\infty)$, $K_0>0$, and $[w]_{A_p}\leq K_0$. Then for any $u\in C_c^\infty(\bR^d)$ and $\varepsilon>0$,
\begin{equation} \label{eq10310054}
\|(-\Delta)^{\sigma'/2}u\|_{L_{p,w}} \leq \varepsilon\|(-\Delta)^{\sigma/2}u\|_{L_{p,w}} + N\varepsilon^{-\sigma'/(\sigma-\sigma')}\|u\|_{L_{p,w}},
\end{equation}
where $N=N(d,p,K_0,\sigma,\sigma')$. 
Moreover, $N$ can be chosen to depend only on $d,p,K_0,\widetilde{\sigma},\sigma'$ when $\sigma=1$, or when $\sigma'<\sigma_0$ and $\sigma$ belongs to either of the other two ranges in \eqref{eq_sigma}.
\end{lemma}

\begin{proof}
For $c>0$, we define the Fourier multiplier operator $T_c$ by
\begin{equation*}
T_cu(x) := \cF^{-1}\left[m_c\hat{u}\right](x),
\end{equation*}
where
\begin{equation*}
m_c(\xi) := \frac{c^{\sigma-\sigma'}|\xi|^{\sigma'}}{c^\sigma+|\xi|^\sigma}.
\end{equation*}
Then one can verify that for every multi-index $\alpha$ with $|\alpha|\leq d+2$,
\begin{equation} \label{eq11082252}
\sup_{\xi\neq0} |\xi|^{|\alpha|} \left|D^\alpha m_c(\xi)\right| \leq N(d),
\end{equation}
uniformly for $c>0$ and $0<\sigma'<\sigma<2$. 
Thus, by the weighted Mikhlin multiplier theorem (see \cite{KW79}),
\begin{equation*}
\|T_c u\|_{L_{p,w}}\leq N\|u\|_{L_{p,w}},
\end{equation*}
where $N=N(d,p,K_0)$. Since
\begin{equation*}
(-\Delta)^{\sigma'/2}u = T_c \left(c^{\sigma'}u + c^{\sigma'-\sigma}(-\Delta)^{\sigma/2}u\right),
\end{equation*}
we obtain
\begin{equation*}
\|(-\Delta)^{\sigma'/2}u\|_{L_{p,w}} \leq Nc^{-(\sigma-\sigma')}\|(-\Delta)^{\sigma/2}u\|_{L_{p,w}} + N c^{\sigma'}\|u\|_{L_{p,w}}.
\end{equation*}

By considering
\begin{equation*}
c = \left(\frac{N}{\varepsilon}\right)^{1/(\sigma-\sigma')},
\end{equation*}
we see that
\begin{equation*}
\|(-\Delta)^{\sigma'/2}u\|_{L_{p,w}} \leq \varepsilon \|(-\Delta)^{\sigma/2}u\|_{L_{p,w}} + N^{\sigma/(\sigma-\sigma')} \varepsilon^{-\sigma'/(\sigma-\sigma')} \|u\|_{L_{p,w}},
\end{equation*}
which proves \eqref{eq10310054}.
Finally, if $\sigma'<\sigma_0\leq\sigma$, then
\begin{equation*}
N^{\sigma/(\sigma-\sigma')} \leq N^{\sigma_0/(\sigma_0-\sigma')}.
\end{equation*}
Here, we may assume that $N\geq1$ without loss of generality.
This yields the asserted dependence of the constants.
The lemma is proved.
\end{proof}

\begin{corollary}
    Let $\sigma\in(1,2)$, $p\in(1,\infty)$, $K_0>0$, and $[w]_{A_p}\leq K_0$. Then for any $u\in C_c^\infty(\bR^d)$ and $\varepsilon>0$,
\begin{equation} \label{eq_grad}
\varepsilon^{1-\sigma} \|D_xu\|_{L_{p,w}} \leq N \|(-\Delta)^{\sigma/2}u\|_{L_{p,w}} + N \varepsilon^{-\sigma} \|u\|_{L_{p,w}},
\end{equation}
where $N=N(d,p,K_0,\widetilde{\sigma})$ and $\widetilde{\sigma}$ is given by \eqref{eq_sigma}.
\end{corollary}

\begin{proof}
By the weighted boundedness of the Riesz transforms (see, e.g. \cite{P08}),
\begin{equation} \label{eq11050057}
\|D_xu\|_{L_{p,w}}\leq N(d,p,K_0)\|(-\Delta)^{1/2}u\|_{L_{p,w}}.
\end{equation}
By \eqref{eq10310054} with $\sigma'=1$ and $\varepsilon$ replaced by $\varepsilon^{\sigma-1}$,
\begin{equation*}
\|(-\Delta)^{1/2}u\|_{L_{p,w}} \leq \varepsilon^{\sigma-1}\|(-\Delta)^{\sigma/2}u\|_{L_{p,w}} + N\varepsilon^{-1}\|u\|_{L_{p,w}}.
\end{equation*}
Multiplying by $\varepsilon^{1-\sigma}$ and applying \eqref{eq11050057}, we obtain \eqref{eq_grad}.
The corollary is proved.
\end{proof}


We next recall the boundedness result in \cite[Proposition 4.1]{DJK23}, together with its version for truncated operators.
For a kernel $K(x,y)$ and $R>0$, we define the truncated operator
\begin{equation} \label{eq_trun}
L_{K,R}u(x) := \int_{|y|<R} \delta_yu(x) K(x,y) \frac{dy}{|y|^{d+\sigma}}.
\end{equation}

\begin{lemma} \label{lem_bound}
Let $\sigma\in(0,2)$, $p\in(1,\infty)$, $K_0>0$, and $[w]_{A_p}\leq K_0$.
Suppose that $K=K(y)$ is a measurable function satisfying Assumption \ref{as_bound}. Then for $u\in H_{p,w}^{\sigma}$ and $R>0$,
\begin{equation} \label{eq10311722}
\|L_Ku\|_{L_{p,w}} + \|L_{K,R}u\|_{L_{p,w}} \leq N\Lambda \|(-\Delta)^{\sigma/2}u\|_{L_{p,w}},
\end{equation}
where $N=N(d,p,K_0,\widetilde{\sigma})$ and $\widetilde{\sigma}$ is defined in \eqref{eq_sigma}.
\end{lemma}

\begin{proof}
See \cite[Proposition 4.1]{DJK23}. We only remark that the kernel of $L_{K,R}$ still satisfies Assumption \ref{as_bound}. The lemma is proved.
\end{proof}

\begin{lemma} \label{lem_simple}
Let $\lambda>0$, $\sigma\in(0,2)$, $p\in(1,\infty)$, $K_0>0$, and $[w]_{A_p}\leq K_0$.
Suppose that $K=K(y)$ is a measurable function satisfying Assumptions \ref{as_bound} and \ref{as_nonde}.
Then for $u\in H_{p,w}^\sigma$ satisfying
\begin{equation} \label{eq11011228}
    L_Ku-\lambda u = f,
\end{equation}
we have
\begin{equation} \label{eq10311738}
\|(-\Delta)^{\sigma/2}u\|_{L_{p,w}} + \lambda \|u\|_{L_{p,w}} \leq N \|f\|_{L_{p,w}},
\end{equation}
where $N=N(d,p,K_0,\nu,\Lambda,\widetilde{\sigma})$ and $\widetilde{\sigma}$ is defined in \eqref{eq_sigma}.
Moreover, for $f\in L_{p,w}$ and $\lambda>0$,  there exists a unique solution $u\in H_{p,w}^{\sigma}$ to \eqref{eq11011228}.
\end{lemma}

\begin{proof}
The estimate \eqref{eq10311738} follows from arguments in the existing literature. We provide a brief explanation.

First, we discuss the following parabolic a priori estimate for $u\in C_c^\infty([0,\infty)\times\bR^d)$ with $u(0,\cdot)=0$:
\begin{align} \label{eq11031609}
    &\|\partial_t u\|_{L_{p}((0,T)\times \bR^d; w(x)dtdx)} + \|(-\Delta)^{\sigma/2}u\|_{L_{p}((0,T)\times \bR^d; w(x)dtdx)} \nonumber
    \\
    &+ \lambda \|u\|_{L_{p}((0,T)\times \bR^d; w(x)dtdx)} \leq N \|\partial_tu - L_K u + \lambda u\|_{L_{p}((0,T)\times \bR^d; w(x)dtdx)},
\end{align}
where $N=N(d,p,\widetilde{\sigma},\nu,\Lambda,K_0)$ is independent of $T$.
The estimate \eqref{eq11031609} for $p=2$ and $w=1$ follows by repeating the proof in \cite[Section 3]{DL23}.
The only difference is that our nondegeneracy condition \eqref{eq_nonde} is more general than that imposed in \cite{DL23}. However, the same argument can be applied by replacing inequality (3.8) therein with the following inequality.
Since $1-\cos{t} \geq t^2/3$ for $|t|\leq1$, by \eqref{eq_nonde},
\begin{align} \label{eq11041328}
\int_{\bR^d} \left(1-\cos(\xi\cdot y)\right) K(y) \frac{dy}{|y|^{d+\sigma}}
\geq \frac{1}{3} \int_{|\xi\cdot y|\leq1} |\xi\cdot y|^2 K(y) \frac{dy}{|y|^{d+\sigma}} \geq \frac{\nu}{3}|\xi|^\sigma.
\end{align}
For general $p$ and $w$, we can follow the arguments in Sections 3 and 4 of \cite{DL23} word for word. It is worth noting that, in the derivation of \eqref{eq11031609}, the lower bound \eqref{eq_nonde} is used only to obtain \eqref{eq11041328}. For the dependence of the constant, we refer to \cite[Remark 2.11]{DL23}.

Next, we prove \eqref{eq10311738}. 
Due to Remark \ref{rem_dense} and Lemma \ref{lem_bound}, we may assume that $u\in C_c^\infty(\bR^d)$. Let $f:=L_Ku-\lambda u$.
Let $\eta\in C_c^\infty((0,1))$ be a nonzero function. For $n>0$, $u_n(t,x):=\eta(t/n)u(x)$ satisfies
\begin{equation*}
\partial_t u_n-L_Ku_n + \lambda u_n = n^{-1}\eta'(t/n)u - \eta(t/n) f.
\end{equation*}
By applying \eqref{eq11031609} with $T=n$ to $u_n$ and dividing by $n^{1/p}\|\eta\|_{L_p(0,1)}$, we obtain
\begin{equation*}
\|(-\Delta)^{\sigma/2}u\|_{L_{p,w}} + \lambda\|u\|_{L_{p,w}} \leq N \|f\|_{L_{p,w}} + \frac{N}{n}\|u\|_{L_{p,w}}.
\end{equation*}
Thus, \eqref{eq10311738} follows by letting $n\to\infty$.


Lastly, we prove the solvability.
Since the uniqueness follows from \eqref{eq10311738}, we only need to prove the existence.
Due to \eqref{eq10311738} and the method of continuity, it suffices to consider the case $L_K=-(-\Delta)^{\sigma/2}$. Moreover, due to the denseness of $C_c^\infty(\bR^d)$ in $L_{p,w}$ and Lemma \ref{lem_bound}, we may assume that $f\in C_c^\infty(\bR^d)$.

Define
\begin{equation*}
u := -\cF^{-1}\left[\frac{\hat{f}}{\lambda+|\cdot|^\sigma}\right].
\end{equation*}
Since $\hat{f}\in\cS(\bR^d)$ and $\lambda>0$, we have $(1+|\cdot|)^k\hat{u}\in L_1(\bR^d)\cap L_2(\bR^d)$ for any $k\geq0$.
Thus, $u$ is infinitely differentiable, $u$ and its derivatives are bounded and belong to $L_2(\bR^d)$, and
\begin{equation*}
-(-\Delta)^{\sigma/2}u-\lambda u = -\cF^{-1}\left[(|\cdot|^\sigma+\lambda)\hat{u}\right] = f.
\end{equation*}
Therefore, it remains to show that $u\in H_{p,w}^\sigma$.
Let us consider the Fourier multiplier operator
\begin{equation*}
Tf := \cF^{-1} \left[m\hat{f}\right],
\end{equation*}
where $m(\xi) := \frac{(1+|\xi|^2)^{\sigma/2}}{\lambda+|\xi|^\sigma}$.
Then $(1-\Delta)^{\sigma/2}u=-Tf$.
Here, one can verify that $m$ satisfies the Mikhlin condition \eqref{eq11082252} with a constant $N$ possibly depending on $d$ and $\lambda$.
Hence, by the weighted Mikhlin multiplier theorem,
\begin{equation*}
\|u\|_{H_{p,w}^\sigma} = \|Tf\|_{L_{p,w}} \leq N\|f\|_{L_{p,w}},
\end{equation*}
where $N=N(d,p,K_0,\lambda)$.
In other words, $u\in H_{p,w}^\sigma$ is a solution to \eqref{eq11011228} with $L_K=-(-\Delta)^{\sigma/2}$.
The lemma is proved.
\end{proof}

\section{Analysis of nonlocal operators} \label{sec_analysis}

In this section, we develop some useful estimates and techniques for the analysis of nonlocal operators.

We begin with the following uniform estimate for $\delta_yu$ over dyadic annuli.
Recall from \eqref{eq_differ} that
\begin{equation} \label{eq11071440}
    \delta_yu(x) := u(x+y)-u(x)-y^{(\sigma)}\cdot\nabla u(x).
\end{equation}

\begin{lemma}
Let $r>0$, $\sigma\in(0,2)$, $p\in(1,\infty)$, $K_0>0$, and $[w]_{A_p}\leq K_0$.
For $u\in C_c^\infty(\bR^d)$,
\begin{equation} \label{eq_shell}
(2-\sigma) \left\|\int_{\{r/2\leq|y|<r\}} \frac{|\delta_yu|}{|y|^{d+\sigma}} dy\right\|_{L_{p,w}} \leq N \|(-\Delta)^{\sigma/2}u\|_{L_{p,w}},
\end{equation}
where $N=N(d,p,K_0,\widetilde{\sigma})$ is independent of $r$.
\end{lemma}

\begin{proof}
We divide the proof into cases according to the value of $\sigma$.

\textbf{1.}
In this case, we deal with $\sigma\in(0,1)$.
By Lemma \ref{lem11041430},
\begin{align} \label{eq11042317}
u(x+y)-u(x) &= N_0 \int_{\bR^d} \left( |x+y-z|^{-d+\sigma}-|x-z|^{-d+\sigma} \right) (-\Delta)^{\sigma/2}u(z) dz \nonumber
\\
&=N_0\int_{|x-z|\leq2|y|} \left(|x+y-z|^{\sigma-d}-|x-z|^{\sigma-d}\right) (-\Delta)^{\sigma/2}u(z) dz \nonumber
\\
&\quad+N_0\int_{|x-z|>2|y|}\left(|x+y-z|^{\sigma-d}-|x-z|^{\sigma-d}\right) (-\Delta)^{\sigma/2}u(z) dz \nonumber
\\
&=: N_0(I_1+I_2),
\end{align}
where $N_0 := \frac{\Gamma((d-\sigma)/2)}{2^\sigma\pi^{d/2}\Gamma(\sigma/2)}$.

We first estimate $I_1$. Since $|x+y-z|\leq |x-z|+|y|\leq3|y|$ if $|x-z|\leq 2|y|$, by \eqref{eq10291653},
\begin{align} \label{eq11042315}
|I_1| &\leq \int_{|x-z|\leq2|y|} \frac{|(-\Delta)^{\sigma/2}u(z)|}{|x-z|^{d-\sigma}} dz + \int_{|x-z|\leq2|y|} \frac{|(-\Delta)^{\sigma/2}u(z)|}{|x+y-z|^{d-\sigma}} dz \nonumber
\\
&\leq \int_{|x-z|\leq2|y|} \frac{|(-\Delta)^{\sigma/2}u(z)|}{|x-z|^{d-\sigma}} dz + \int_{|x+y-z|\leq3|y|} \frac{|(-\Delta)^{\sigma/2}u(z)|}{|x+y-z|^{d-\sigma}} dz \nonumber
\\
&\leq \frac{N}{\sigma} \left((2|y|)^\sigma \bM((-\Delta)^{\sigma/2}u)(x)+(3|y|)^\sigma \bM((-\Delta)^{\sigma/2}u)(x+y)\right) \nonumber
\\
&\leq N(\sigma_0) |y|^\sigma \left(\bM ((-\Delta)^{\sigma/2}u)(x) + \bM ((-\Delta)^{\sigma/2}u)(x+y)\right),
\end{align}
where $\sigma_0$ is introduced in \eqref{eq_sigma}.

Next, we deal with $I_2$. For $|x-z|>2|y|$ and $t\in[0,1]$,
\begin{equation*}
|x+ty-z|\geq|x-z|-t|y|\geq|x-z|-|y|\geq\frac{1}{2}|x-z|.
\end{equation*}
Therefore, by the mean value theorem,
\begin{equation*}
\left||x+y-z|^{\sigma-d}-|x-z|^{\sigma-d}\right| \leq N(d-\sigma)|y||x-z|^{\sigma-d-1}.
\end{equation*}
This together with \eqref{eq10291650} yields
\begin{equation} \label{eq11042316}
    |I_2| \leq N \frac{d-\sigma}{1-\sigma} |y|^\sigma \bM((-\Delta)^{\sigma/2}u)(x).
\end{equation}
Thus, by \eqref{eq11042317}, \eqref{eq11042315}, and \eqref{eq11042316},
\begin{align} \label{eq11050019}
&|u(x+y)-u(x)| \nonumber
\\
&\leq N_0 \left(N(\sigma_0) + N\frac{d-\sigma}{1-\sigma}\right) |y|^\sigma \left(\bM((-\Delta)^{\sigma/2}u)(x) + \bM((-\Delta)^{\sigma/2}u)(x+y)\right),
\end{align}
which yields
\begin{equation} \label{eq11042335}
\int_{\{r/2\leq|y|<r\}} \frac{|\delta_yu(x)|}{|y|^{d+\sigma}}dy \leq N \left(\bM((-\Delta)^{\sigma/2}u)(x)+\bM(\bM((-\Delta)^{\sigma/2}u))(x)\right).
\end{equation}
Here, we used
\begin{align} \label{eq11042324}
\int_{\{r/2\leq|y|<r\}} \frac{|\bM((-\Delta)^{\sigma/2}u)(x+y)|}{|y|^d} dy &\leq \left(\frac{2}{r}\right)^d\int_{B_r(x)}|\bM((-\Delta)^{\sigma/2}u)(z)| dz \nonumber
\\
&\leq N \bM(\bM ((-\Delta)^{\sigma/2}u))(x).
\end{align}
By applying the weighted maximal function theorem (see \cite[Theorem 2.2]{DK18}) to \eqref{eq11042335}, we obtain \eqref{eq_shell}.

\textbf{2.}
Next, we assume that $\sigma\in(1,2)$.
Let
\begin{equation*}
v := (-\Delta)^{\alpha/2}D_xu,
\end{equation*}
where $\alpha:=\sigma-1\in(0,1)$.
Since $D_xu\in C_c^\infty(\bR^d)$, by Lemma \ref{lem11041430} with $\sigma=\alpha$,
\begin{equation*}
D_xu(x) = \frac{\Gamma((d+1-\sigma)/2)}{2^{\sigma-1}\pi^{d/2} \Gamma((\sigma-1)/2)} \int_{\bR^d} |x-z|^{\alpha-d} v(z) dz.
\end{equation*}
By applying \eqref{eq11050019} to $D_xu$,
\begin{align*}
|D_xu(x+y)-D_xu(x)| \leq \tilde{N} |y|^\alpha \left(\bM(|v|)(x)+\bM(|v|)(x+y)\right),
\end{align*}
where
\begin{equation*}
    \tilde{N}:= N \frac{\Gamma((d+1-\sigma)/2)}{2^{\sigma-1}\pi^{d/2} \Gamma((\sigma-1)/2)} \left( N(\sigma_0) + \frac{d+1-\sigma}{2-\sigma} \right).
\end{equation*}
Here, one can find that $(2-\sigma)\tilde{N}\leq N(\widetilde{\sigma})$, where $\widetilde{\sigma}$ is given by \eqref{eq_sigma}.
Thus,
\begin{equation} \label{eq11050051}
    (2-\sigma)|D_xu(x+y)-D_xu(x)| \leq N(\widetilde{\sigma}) |y|^\alpha \left(\bM(|v|)(x) + \bM(|v|)(x+y)\right).
\end{equation}

By the fundamental theorem of calculus,
\begin{equation*}
\delta_yu(x)=\int_0^1\left(D_xu(x+ty)-D_xu(x)\right)\cdot ydt.
\end{equation*}
This and \eqref{eq11050051} yield
\begin{align*}
    &(2-\sigma)\int_{\{r/2\leq |y|<r\}} \frac{|\delta_yu(x)|}{|y|^{d+\sigma}} dy
\\
&\leq N \int_0^1 t^{\sigma-1} \int_{\{r/2\leq |y|<r\}} \frac{\bM(|v|)(x+ty)+\bM(|v|)(x)}{|y|^d} dydt.
\end{align*}
For $t>0$, as in \eqref{eq11042324},
\begin{align*}
\int_{\{r/2\leq |y|<r\}} \frac{\bM(|v|)(x+ty)}{|y|^d}dy &= \int_{\{tr/2\leq |z|<tr\}} \frac{\bM(|v|)(x+z)}{|z|^d} dz \leq N \bM(\bM(|v|))(x).
\end{align*}
Hence, we get
\begin{equation} \label{eq11050110}
(2-\sigma)\int_{\{r/2\leq |y|<r\}}\frac{|\delta_yu(x)|}{|y|^{d+\sigma}}dy \leq N \left(\bM(|v|)(x)+\bM(\bM(|v|))(x)\right).
\end{equation}
By the weighted maximal function theorem,
\begin{equation*}
    (2-\sigma) \left\|\int_{\{r/2\leq|y|<r\}} \frac{|\delta_yu|}{|y|^{d+\sigma}} dy\right\|_{L_{p,w}} \leq N \|(-\Delta)^{\alpha/2} D_xu\|_{L_{p,w}}.
\end{equation*}
Now \eqref{eq_shell} follows from a weighted result for the Riesz transforms \eqref{eq11050057} and the identity $(-\Delta)^{\sigma/2}=(-\Delta)^{1/2}(-\Delta)^{\alpha/2}$.

 \textbf{3.}
 Lastly, we consider the case when $\sigma=1$.
 By the fundamental theorem of calculus,
\begin{equation*}
|\delta_yu(x)| \leq |y| \int_0^1|D_xu(x+ty)| dt + |y|1_{\{|y|<1\}} |D_xu(x)|.
\end{equation*}
As in \eqref{eq11050110},
\begin{equation} \label{eq11050111}
\int_{\{r/2\leq |y|<r\}} \frac{|\delta_yu(x)|}{|y|^{d+1}}dy\leq N \left( \bM(\bM(|D_xu|))(x)+|D_xu(x)|\right).
\end{equation}
By the weighted maximal function theorem and \eqref{eq11050057},
\begin{equation*}
\|\bM(\bM(|D_xu|))\|_{L_{p,w}} + \|D_xu\|_{L_{p,w}} \leq N \|(-\Delta)^{1/2}u\|_{L_{p,w}}.
\end{equation*}
This and \eqref{eq11050111} yield \eqref{eq_shell}.
The lemma is proved.
\end{proof}

In order to carry out our analysis on a spatial grid, we introduce the following notation.
For $r\in(0,1]$ and $m\in\bZ^d$, let
\begin{equation*}
x_{r,m}:=\frac{r}{\sqrt d}m.
\end{equation*}
Let $\eta\in C_c^\infty(B_1)$ and $\widetilde{\eta}\in C_c^\infty(B_4)$ such that
\begin{equation} \label{eq11090037}
    0\leq\eta,\widetilde{\eta}\leq1, \quad \widetilde{\eta}=1 \text{ on } B_2,
\end{equation}
and
\begin{equation*}
\sum_{m\in\bZ^d} \eta \left(x-\frac{m}{\sqrt d}\right)=1, \quad x\in\bR^d.
\end{equation*}
We note that these functions can be chosen depending only on $d$. Define
\begin{equation} \label{eqeta}
\eta_{r,m}(x) := \eta \left(\frac{x-x_{r,m}}{r}\right), \quad \widetilde{\eta}_{r,m}(x) := \widetilde{\eta} \left(\frac{x-x_{r,m}}{2r}\right).
\end{equation}
Then
\begin{align} \label{eq11091424}
&0\leq\eta_{r,m}\leq1, \quad \sum_{m\in\bZ^d}\eta_{r,m}=1, \quad \text{supp}(\eta_{r,m}) \subset B_r(x_{r,m}), \nonumber
\\
&\widetilde{\eta}_{r,m}=1 \text{ on } B_{4r}(x_{r,m}), \quad \text{supp}(\widetilde{\eta}_{r,m}) \subset B_{8r}(x_{r,m}),
\end{align}
and, for $k\geq0$,
\begin{equation*}
|D_x^k\eta_{r,m}|+|D_x^k\widetilde{\eta}_{r,m}| \leq N(k)r^{-k}.
\end{equation*}

\begin{lemma} \label{lem11141719}
Let $r\in(0,1]$, $p\in(1,\infty)$, and let $w$ be a nonnegative locally integrable weight on $\bR^d$. Then for any function $u$,
\begin{equation} \label{eq11061351}
N(d,p) \|u\|_{L_{p,w}} \leq \left(\sum_{m\in\bZ^d} \left\|\eta_{r,m}u\right\|_{L_{p,w}}^p\right)^{1/p} \leq \|u\|_{L_{p,w}}.
\end{equation}
Moreover,
\begin{equation} \label{eq11061352}
\sum_{m\in\bZ^d} \left\|\widetilde{\eta}_{r,m}u\right\|_{L_{p,w}}^p \leq N(d) \|u\|_{L_{p,w}}^p.
\end{equation}
\end{lemma}

\begin{proof}
Let $\kappa\geq1$ and $m_0,m\in \bZ^d$. Suppose that $B_{\kappa r}(x_{r,m}) \cap B_{\kappa r}(x_{r,m_0}) \neq \emptyset$ and take $x\in B_{\kappa r}(x_{r,m}) \cap B_{\kappa r}(x_{r,m_0})$. Then we have
\begin{align*}
|m_0-m| = \frac{\sqrt{d}}{r}|x_{r,m_0}-x_{r,m}| \leq \frac{\sqrt{d}}{r} \left(|x_{r,m_0}-x|+|x-x_{r,m}|\right) < 2\kappa \sqrt{d}.
\end{align*}
Thus, one can deduce that
\begin{equation} \label{eq11061527}
\sum_{m \in \bZ^d} 1_{B_{\kappa r}(x_{r,m})}(x) \leq (4\kappa\sqrt{d}+1)^d.
\end{equation}

Now we prove \eqref{eq11061351}. 
Since $\text{supp}(\eta_{r,m}) \subset B_r(x_{r,m})$, by H\"older's inequality and \eqref{eq11061527} with $\kappa=1$,
\begin{align*}
1=\sum_m\eta_{r,m}(x) &= \sum_{m\in\bZ^d} \eta_{r,m}(x) 1_{B_r(x_{r,m})}(x)
\\
&\leq\left(\sum_{m\in\bZ^d} \eta_{r,m}(x)^p\right)^{1/p}\left(\sum_{m\in\bZ^d} 1_{B_r(x_{r,m})}(x)\right)^{1-1/p}
\\
&\leq N(d,p) \left(\sum_{m\in\bZ^d} \eta_{r,m}(x)^p\right)^{1/p},
\end{align*}
where $N(d,p)>0$.
On the other hand, due to $0\leq\eta_{r,m}\leq1$,
\begin{equation*}
\sum_{m\in\bZ^d} \eta_{r,m}(x)^p \leq \sum_{m\in\bZ^d} \eta_{r,m}(x) = 1.
\end{equation*}
Thus, by Fubini's theorem, we get \eqref{eq11061351}.

For \eqref{eq11061352}, by \eqref{eq11061527} with $\kappa=8$,
\begin{equation*}
\sum_{m\in \bZ^d} \widetilde{\eta}_{r,m}(x)^p \leq \sum_{m\in\bZ^d} 1_{B_{8r}(x_{r,m})}(x) \leq N(d).
\end{equation*}
Thus, combining this with Fubini's theorem yields \eqref{eq11061352}.
The lemma is proved.
\end{proof}

We next derive an estimate for the commutator terms.

\begin{lemma} \label{lem11131905}
Let $r\in(0,1]$, $\sigma\in(0,2)$, $p\in(1,\infty)$, $K_0>0$, $[w]_{A_p}\leq K_0$, $\kappa\geq1$, and $\zeta_m\in C_c^\infty(B_{\kappa r}(x_{r,m}))$ satisfy $0\leq\zeta_m\leq1$ and $|D_x^k\zeta_m| \leq Nr^{-k}$ for $k=1,2$. Let $K_m=K_m(x,y)$ be measurable kernels satisfying Assumption \ref{as_bound} with the same constant $\Lambda$. Then for $u\in C_c^\infty(\bR^d)$ and $\varepsilon\in(0,1)$,
\begin{equation} \label{eq11071436}
\left(\sum_{m\in\bZ^d} \left\|[L_{K_m},\zeta_m]u\right\|_{L_{p,w}}^p\right)^{1/p} \leq N \Lambda \left( \varepsilon \|(-\Delta)^{\sigma/2}u\|_{L_{p,w}} + N_\varepsilon r^{-\sigma} \|u\|_{L_{p,w}}\right),
\end{equation}
where
\begin{equation} \label{eq11081725}
[L_{K_m},\zeta_m]u := L_{K_m}(\zeta_mu)-\zeta_mL_{K_m}u.
\end{equation}
Here, the constants depend only on $d,p,K_0,\widetilde{\sigma},\kappa$, with an additional dependence on $\varepsilon$ for $N_\varepsilon$, where $\widetilde{\sigma}$ is given by \eqref{eq_sigma}.
\end{lemma}

\begin{proof}
By replacing $K_m$ by $K_m/\Lambda$, we may assume that $\Lambda=1$.
Due to \eqref{intro_op} and \eqref{eq11071440},
\begin{align*}
&[L_{K_m},\zeta_m]u(x)
\\
&= \int_{\bR^d} \left[\left(\zeta_m(x+y)-\zeta_m(x)\right)u(x+y) -y^{(\sigma)}\cdot D_x\zeta_m(x)u(x)\right] K_m(x,y) \frac{dy}{|y|^{d+\sigma}}.
\end{align*}
If $x\notin B_{2\kappa r}(x_{r,m})$ and $|y|<r$, then $|x+y-x_{r,m}|\geq |x-x_{r,m}|-|y|>2\kappa r-r\geq\kappa r$, where we used $\kappa\geq1$. Since $\text{supp}(\zeta_m)\subset B_{\kappa r}(x_{r,m})$, we get $\zeta_m(x+y)=\zeta_m(x)=D_x\zeta_m(x)=0$, which leads to
\begin{equation*}
1_{B_{2\kappa r}(x_{r,m})^c}(x) \int_{|y|<r}\cdots dy = 0.
\end{equation*}
Therefore, we see that
\begin{align} \label{eq11071441}
&[L_{K_m},\zeta_m]u(x) \nonumber
\\
&= \int_{\bR^d} \left[\left(\zeta_m(x+y)-\zeta_m(x)\right)u(x+y) -y^{(\sigma)}\cdot D_x\zeta_m(x)u(x)\right] K_m(x,y) \frac{dy}{|y|^{d+\sigma}} \nonumber
\\
&= 1_{B_{2\kappa r}(x_{r,m})^c}(x) \int_{|y|\geq r} \cdots dy + 1_{B_{2\kappa r}(x_{r,m})}(x) \int_{|y|\geq r} \cdots dy \nonumber
\\
&\quad+ 1_{B_{2\kappa r}(x_{r,m})}(x) \int_{|y|<r} \cdots dy \nonumber
\\
&=: I_{1,m}(x) + I_{2,m}(x) + I_{3,m}(x).
\end{align}

\textbf{Step 1.} We estimate $I_{1,m}$.
Note that if $x\in B_{2\kappa r}(x_{r,m})^c$ and $\zeta_m(x+y)\neq0$, then $x+y\in B_{\kappa r}(x_{r,m})$ and
\begin{equation*}
|y|\geq|x-x_{r,m}|-\kappa r\geq\frac{|x-x_{r,m}|}{2}\geq\kappa r\geq r.
\end{equation*}
Thus, by the change of variables $z=x+y$,
\begin{align} \label{eq11071556}
|I_{1,m}(x)| &= \left| 1_{B_{2\kappa r}(x_{r,m})^c}(x) \int_{|y|\geq r}\zeta_m(x+y) u(x+y) K_m(x,y) \frac{dy}{|y|^{d+\sigma}} \right| \nonumber
\\
&\leq N(2-\sigma) |x-x_{r,m}|^{-d-\sigma} \int_{B_{\kappa r}(x_{r,m})} |u(z)| dz \nonumber
\\
&\leq N(2-\sigma) |x-x_{r,m}|^{-d-\sigma} (\kappa r)^d \inf_{z\in B_{\kappa r}(x_{r,m})}\bM u(z).
\end{align}

For a measurable set $E$, we write $w(E):=\int_Ew(x)dx$. By H\"older's inequality and the $A_p$ condition, for any ball $B\subset\bR^d$ and $c\geq1$,
\begin{equation*}
|B|^p \leq w(B) \left(\int_{cB}w^{-1/(p-1)} dx \right)^{p-1}\leq K_0 |cB|^p \frac{w(B)}{w(cB)},
\end{equation*}
which yields
\begin{equation*}
w(cB) \leq K_0 c^{dp} w(B).
\end{equation*}
This leads to
\begin{align*}
&\int_{B_{2\kappa r}(x_{r,m})^c}|x-x_{r,m}|^{-(d+\sigma)p} w(x) dx
\\
&= \sum_{k\geq1} \int_{2^k\kappa r\leq|x-x_{r,m}|<2^{k+1}\kappa r} |x-x_{r,m}|^{-(d+\sigma)p} w(x) dx
\\
&\leq\sum_{k\geq1} (2^k\kappa r)^{-(d+\sigma)p} w(B_{2^{k+1}\kappa r}(x_{r,m}))
\\
&\leq K_0 2^{dp} (\kappa r)^{-(d+\sigma)p} w(B_{\kappa r}(x_{r,m})) \sum_{k\geq1}2^{-k\sigma p}
\\
&\leq N(K_0,d,p,\widetilde{\sigma}) (\kappa r)^{-(d+\sigma)p} w(B_{\kappa r}(x_{r,m})).
\end{align*}
This and \eqref{eq11071556} yield
\begin{align*}
\|I_{1,m}\|_{L_{p,w}}^p &\leq N r^{-\sigma p} w(B_{\kappa r}(x_{r,m})) \left(\inf_{B_{\kappa r}(x_{r,m})}\bM u\right)^p
\\
&\leq N r^{-\sigma p} \int_{B_{\kappa r}(x_{r,m})}(\bM u(x))^p w(x) dx.
\end{align*}
Hence, by \eqref{eq11061527} and the weighted maximal function theorem, we obtain
\begin{equation} \label{eq11071357}
\sum_{m\in\bZ^d}\|I_{1,m}\|_{L_{p,w}}^p \leq N r^{-\sigma p} \|u\|_{L_{p,w}}^p.
\end{equation}

\textbf{Step 2.} In this step, we consider $I_{2,m}$.
Due to $0\leq\zeta_m\leq1$ and \eqref{eq10291650},
\begin{align*}
&\int_{|y|\geq r} \left|\zeta_m(x+y)-\zeta_m(x)\right||u(x+y)|K_m(x,y)\frac{dy}{|y|^{d+\sigma}}
\\
&\leq 2(2-\sigma) \int_{|y|\geq r}\frac{|u(x+y)|}{|y|^{d+\sigma}} dy \leq N \frac{2-\sigma}{\sigma} r^{-\sigma} \bM u(x).
\end{align*}

We next estimate the terms in $I_{2,m}$ containing $y^{(\sigma)}\cdot D_x\zeta_m(x)u(x)$. Since $y^{(\sigma)}=0$ when $\sigma<1$, we only consider the case when $\sigma\geq1$.
If $\sigma=1$, then \eqref{eq_cancel} leads to
\begin{equation*}
u(x)D_x\zeta_m(x) \cdot \int_{r\leq|y|<1} y K_m(x,y)\frac{dy}{|y|^{d+1}}=0.
\end{equation*}
 If $\sigma>1$, then by $|D_x\zeta_m|\leq N/r$,
\begin{align*}
\left|\int_{|y|\geq r} y\cdot D_x\zeta_m(x) u(x) K_m(x,y) \frac{dy}{|y|^{d+\sigma}}\right| &\leq \frac{N(2-\sigma)}{r} |u(x)| \int_{|y|\geq r} |y|^{1-d-\sigma}dy
\\
&\leq N\frac{2-\sigma}{\sigma-1}r^{-\sigma}|u(x)|.
\end{align*}
 Consequently,
\begin{equation*}
|I_{2,m}(x)|\leq N r^{-\sigma} 1_{B_{2\kappa r}(x_{r,m})}(x) \left(\bM u(x)+|u(x)|\right).
\end{equation*}
By \eqref{eqmaximal}, and \eqref{eq11061527} with $2\kappa$ in place of $\kappa$,
\begin{equation} \label{eq11081147}
\left(\sum_{m\in\bZ^d} \|I_{2,m}\|_{L_{p,w}}^p\right)^{1/p}\leq Nr^{-\sigma}\|u\|_{L_{p,w}}.
\end{equation}

\textbf{Step 3-1.} We consider $I_{3,m}$ when $\sigma\in(0,1)$.
Using $\left|\zeta_m(x+y)-\zeta_m(x)\right|\leq N|y|/r$ and \eqref{eq10291653} with $\gamma=1-\sigma$,
\begin{align*}
|I_{3,m}(x)| &\leq \frac{N(2-\sigma)}{r}1_{B_{2\kappa r}(x_{r,m})}(x) \int_{|y|<r} |u(x+y)| \frac{dy}{|y|^{d-(1-\sigma)}}
\\
&\leq N \frac{2-\sigma}{1-\sigma}r^{-\sigma} 1_{B_{2\kappa r}(x_{r,m})}(x) \bM u(x).
\end{align*}
By \eqref{eqmaximal}, and \eqref{eq11061527} with $2\kappa$ in place of $\kappa$,
\begin{equation} \label{eq11081148}
\left(\sum_{m\in\bZ^d} \|I_{3,m}\|_{L_{p,w}}^p\right)^{1/p} \leq Nr^{-\sigma} \|u\|_{L_{p,w}}.
\end{equation}

\textbf{Step 3-2.} We deal with $I_{3,m}$ when $\sigma\in(1,2)$.
We use the identity
\begin{align} \label{eq11071640}
&\left(\zeta_m(x+y)-\zeta_m(x)\right)u(x+y)-y\cdot D_x\zeta_m(x)u(x) \nonumber
\\
&= \left(\zeta_m(x+y)-\zeta_m(x)-y\cdot D_x\zeta_m(x)\right) u(x) \nonumber
\\
&\quad+ \left(\zeta_m(x+y)-\zeta_m(x)\right) \left(u(x+y)-u(x)\right).
\end{align}
Using $|D_x^2 \zeta_m|\leq Nr^{-2}$,
\begin{align} \label{eq11071641}
&\int_{|y|<r} \left|\zeta_m(x+y)-\zeta_m(x)-y\cdot D_x\zeta_m(x)\right| |u(x)| K_m(x,y) \frac{dy}{|y|^{d+\sigma}} \nonumber
\\
&\leq N(2-\sigma)r^{-2} |u(x)| \int_{|y|<r} |y|^{2-d-\sigma} dy =N r^{-\sigma} |u(x)|.
\end{align}
Next, since $|\zeta_m(x+y)-\zeta_m(x)|\leq\frac{N}{r}|y|$ and
\begin{equation*}
|u(x+y)-u(x)| \leq |y| \int_0^1 |D_xu(x+ty)| dt,
\end{equation*}
we get
\begin{align} \label{eq11071642}
&\int_{|y|<r} \left|\zeta_m(x+y)-\zeta_m(x)\right| \left|u(x+y)-u(x)\right| K_m(x,y) \frac{dy}{|y|^{d+\sigma}} \nonumber
\\
& \leq N \frac{(2-\sigma)}{r}\int_0^1\int_{|y|<r}|D_xu(x+ty)|\frac{dy}{|y|^{d-(2-\sigma)}}dt.
\end{align}
By the change of variables $z=ty$ and \eqref{eq10291653} with $\gamma=2-\sigma$,
\begin{align} \label{eq11071643}
\int_{|y|<r}|D_xu(x+ty)|\frac{dy}{|y|^{d-(2-\sigma)}} &= t^{\sigma-2}\int_{|z|<tr} |D_xu(x+z)| \frac{dz}{|z|^{d-(2-\sigma)}} \nonumber
\\
&\leq\frac{N}{2-\sigma}t^{\sigma-2} (tr)^{2-\sigma} \bM(|D_xu|)(x) \nonumber
\\
&=\frac{N}{2-\sigma} r^{2-\sigma} \bM(|D_xu|)(x).
\end{align}
Combining \eqref{eq11071640}, \eqref{eq11071641}, \eqref{eq11071642}, and \eqref{eq11071643},
\begin{equation*}
|I_{3,m}(x)|\leq N1_{B_{2\kappa r}(x_{r,m})}(x)\left(r^{-\sigma}|u(x)|+r^{1-\sigma}\bM(|D_xu|)(x)\right).
\end{equation*}
Thus, by \eqref{eqmaximal}, and \eqref{eq11061527} with $2\kappa$ in place of $\kappa$,
\begin{equation*}
\left(\sum_{m\in\bZ^d} \left\|I_{3,m}\right\|_{L_{p,w}}^p\right)^{1/p}\leq N\left(r^{-\sigma}\left\|u\right\|_{L_{p,w}}+r^{1-\sigma}\left\|D_xu\right\|_{L_{p,w}}\right).
\end{equation*}
By \eqref{eq_grad} with $\varepsilon^{1/(\sigma-1)} r$ in place of $\varepsilon$, we obtain
\begin{equation*}
r^{1-\sigma} \|D_xu\|_{L_{p,w}} \leq N\varepsilon \|(-\Delta)^{\sigma/2}u\|_{L_{p,w}} + N_{\varepsilon }r^{-\sigma} \|u\|_{L_{p,w}},
\end{equation*}
which yields
\begin{equation} \label{eq11081149}
\left(\sum_{m\in\bZ^d} \left\|I_{3,m}\right\|_{L_{p,w}}^p\right)^{1/p}\leq N\varepsilon \|(-\Delta)^{\sigma/2}u\|_{L_{p,w}} + N_{\varepsilon} r^{-\sigma} \|u\|_{L_{p,w}}.
\end{equation}

\textbf{Step 3-3.} We estimate $I_{3,m}$ when $\sigma=1$.
By using \eqref{eq11071640}, \eqref{eq11071641}, and $|\zeta_m(x+y)-\zeta_m(x)| \leq \frac{N}{r}|y|$,
\begin{equation} \label{eq11081134}
|I_{3,m}(x)|\leq Nr^{-1} 1_{B_{2\kappa r}(x_{r,m})}(x) \left(|u(x)|+J(x)\right),
\end{equation}
where
\begin{equation*}
J(x) := \int_{|y|<r} \frac{|u(x+y)-u(x)|}{|y|^d} dy.
\end{equation*}
By \eqref{eq11050019} with $\sigma=1/2$,
\begin{equation*}
J(x) \leq N \int_{|y|<r}|y|^{1/2-d} \left(\bM((-\Delta)^{1/4}u)(x)+\bM((-\Delta)^{1/4}u)(x+y)\right)dy.
\end{equation*}
Hence, by \eqref{eq10291653} with $\gamma=1/2$,
\begin{equation*}
J(x) \leq Nr^{1/2} \left(\bM((-\Delta)^{1/4}u)(x) + \bM(\bM((-\Delta)^{1/4}u))(x)\right).
\end{equation*}
Now we apply \eqref{eqmaximal} to get
\begin{equation*}
r^{-1} \|J\|_{L_{p,w}} \leq Nr^{-1/2} \|(-\Delta)^{1/4}u\|_{L_{p,w}}.
\end{equation*}
By \eqref{eq10310054} with $\sigma=1$, $\sigma'=1/2$, and $\varepsilon$ replaced by $\varepsilon r^{1/2}$,
\begin{equation*}
r^{-1}\|J\|_{L_{p,w}} \leq N \varepsilon \|(-\Delta)^{1/2}u\|_{L_{p,w}} + N \varepsilon^{-1}r^{-1} \|u\|_{L_{p,w}}.
\end{equation*}
By \eqref{eq11081134} and \eqref{eq11061527} with $2\kappa$ in place of $\kappa$,
\begin{align} \label{eq11081150}
\left(\sum_{m\in\bZ^d} \|I_{3,m}\|_{L_{p,w}}^p\right)^{1/p} &\leq Nr^{-1}\left(\|u\|_{L_{p,w}}+\|J\|_{L_{p,w}}\right) \nonumber
\\
&\leq N\varepsilon \|(-\Delta)^{1/2}u\|_{L_{p,w}} + N \varepsilon^{-1}r^{-1} \|u\|_{L_{p,w}}.
\end{align}

\textbf{Step 4.} Lastly, \eqref{eq11071436} follows by considering \eqref{eq11071441}, \eqref{eq11071357}, \eqref{eq11081147}, \eqref{eq11081148}, \eqref{eq11081149}, and \eqref{eq11081150}.
The lemma is proved.
\end{proof}

As a consequence of Lemmas \ref{lem11141719} and \ref{lem11131905}, we obtain the following localization estimates for the fractional Laplacian.

\begin{corollary}
Let $r\in(0,1]$, $\sigma\in(0,2)$, $p\in(1,\infty)$, $K_0>0$, $[w]_{A_p}\leq K_0$, and $u\in C_c^\infty(\bR^d)$. Then
\begin{equation} \label{eq11081702}
\|(-\Delta)^{\sigma/2}u\|_{L_{p,w}}\leq N\left(\sum_{m\in\bZ^d} \|(-\Delta)^{\sigma/2}(\eta_{r,m}u)\|_{L_{p,w}}^p\right)^{1/p} + Nr^{-\sigma}\|u\|_{L_{p,w}},
\end{equation}
where $N=N(d,p,K_0,\widetilde{\sigma})$ and $\widetilde{\sigma}$ is defined in \eqref{eq_sigma}.
Moreover, let $\kappa\geq1$ and $\zeta_m\in C_c^\infty(B_{\kappa r}(x_{r,m}))$ satisfy $0\leq\zeta_m\leq1$ and $|D_x^k\zeta_m|\leq Nr^{-k}$ for $k=1,2$. Then
\begin{equation} \label{eq11081703}
\left(\sum_{m\in\bZ^d} \|(-\Delta)^{\sigma/2}(\zeta_mu)\|_{L_{p,w}}^p\right)^{1/p} \leq N \|(-\Delta)^{\sigma/2}u\|_{L_{p,w}} + Nr^{-\sigma} \|u\|_{L_{p,w}},
\end{equation}
where $N=N(d,p,K_0,\widetilde{\sigma},\kappa)$.
\end{corollary}

\begin{proof}
To prove \eqref{eq11081702}, we use the identity
\begin{equation*}
\eta_{r,m}(-\Delta)^{\sigma/2}u=(-\Delta)^{\sigma/2}(\eta_{r,m}u)-[(-\Delta)^{\sigma/2},\eta_{r,m}]u,
\end{equation*}
where $[\cdot,\cdot]$ is defined as in \eqref{eq11081725}. By \eqref{eq11061351} and the triangle inequality, we obtain that for $\varepsilon\in(0,1)$,
\begin{align*}
\|(-\Delta)^{\sigma/2}u\|_{L_{p,w}} &\leq N\left(\sum_{m\in\bZ^d}\|\eta_{r,m}(-\Delta)^{\sigma/2}u\|_{L_{p,w}}^p\right)^{1/p}
\\
&\leq N \left(\sum_{m\in\bZ^d} \|(-\Delta)^{\sigma/2}(\eta_{r,m}u)\|_{L_{p,w}}^p\right)^{1/p}
\\
&\quad+ N \left(\sum_{m\in\bZ^d} \|[(-\Delta)^{\sigma/2},\eta_{r,m}]u\|_{L_{p,w}}^p\right)^{1/p}
\\
&\leq N \left(\sum_{m\in\bZ^d} \|(-\Delta)^{\sigma/2}(\eta_{r,m}u)\|_{L_{p,w}}^p\right)^{1/p}
\\
&\quad+ N\varepsilon \|(-\Delta)^{\sigma/2}u\|_{L_{p,w}} + N_\varepsilon r^{-\sigma} \|u\|_{L_{p,w}}.
\end{align*}
Here, for the last inequality, we applied \eqref{eq11071436} with $\zeta_m=\eta_{r,m}$ and $\kappa=1$, and used the fact that
\begin{equation*}
-(-\Delta)^{\sigma/2} = L_{c_{d,\sigma}},
\end{equation*}
where the constant kernel $c_{d,\sigma}=\frac{\sigma2^{\sigma-1}\Gamma((d+\sigma)/2)}{\pi^{d/2}\Gamma(1-\sigma/2)}$ satisfies Assumption \ref{as_bound} with $\Lambda=N(d)$.
By choosing $\varepsilon>0$ sufficiently small so that $N\varepsilon\leq1/2$, we can absorb $\|(-\Delta)^{\sigma/2}u\|_{L_{p,w}}$ into the left-hand side and obtain \eqref{eq11081702}.

Next, we prove \eqref{eq11081703}. Note that
\begin{equation*}
(-\Delta)^{\sigma/2}(\zeta_mu)=\zeta_m(-\Delta)^{\sigma/2}u+[(-\Delta)^{\sigma/2},\zeta_m]u.
\end{equation*}
By \eqref{eq11061527},
\begin{align*}
\sum_{m\in\bZ^d} \|\zeta_m(-\Delta)^{\sigma/2}u\|_{L_{p,w}}^p \leq N \|(-\Delta)^{\sigma/2}u\|_{L_{p,w}}^p.
\end{align*}
Thus, by \eqref{eq11071436} with $\varepsilon=1/2$,
\begin{align*}
\left(\sum_{m\in\bZ^d} \|(-\Delta)^{\sigma/2}(\zeta_mu)\|_{L_{p,w}}^p\right)^{1/p}&\leq\left(\sum_{m\in\bZ^d} \|\zeta_m(-\Delta)^{\sigma/2}u\|_{L_{p,w}}^p\right)^{1/p}
\\
&\quad+\left(\sum_{m\in\bZ^d} \|[(-\Delta)^{\sigma/2},\zeta_m]u\|_{L_{p,w}}^p\right)^{1/p}
\\
&\leq N \|(-\Delta)^{\sigma/2}u\|_{L_{p,w}} + Nr^{-\sigma} \|u\|_{L_{p,w}}.
\end{align*}
Thus, we get \eqref{eq11081703}.
The corollary is proved.
\end{proof}

\begin{lemma} \label{lem11131906}
Let $r\in(0,1]$, $\sigma\in(0,2)$, $p\in(1,\infty)$, $K_0>0$, $[w]_{A_p}\leq K_0$, and $u\in C_c^\infty(\bR^d)$. Let $K_m=K_m(x,y)$ be measurable kernels satisfying Assumption \ref{as_bound} with the same constant $\Lambda$. Then
\begin{equation} \label{eq11081405}
\left( \sum_{m\in\bZ^d} \|\eta_{r,m}L_{K_m}((1-\widetilde{\eta}_{r,m})u)\|_{L_{p,w}}^p \right)^{1/p} \leq N \Lambda r^{-\sigma}\|u\|_{L_{p,w}},
\end{equation}
where $N=N(d,p,K_0,\widetilde{\sigma})$ and $\widetilde{\sigma}$ is defined in \eqref{eq_sigma}.
\end{lemma}

\begin{proof}
Let $v_m := (1-\widetilde{\eta}_{r,m})u$. Then for $x\in\text{supp}(\eta_{r,m})\subset B_r(x_{r,m})$, we see that $\delta_yv_m(x)=v_m(x+y)$, which is zero for $|y|<3r$. Thus, by \eqref{eq_upper} and \eqref{eq10291650},
\begin{align*}
|\eta_{r,m}(x)L_{K_m}v_m(x)| &\leq (2-\sigma)\Lambda \eta_{r,m}(x) \int_{|y|\geq3r} \frac{|u(x+y)|}{|y|^{d+\sigma}} dy
\\
&\leq N \frac{2-\sigma}{\sigma} \Lambda r^{-\sigma} \eta_{r,m}(x) \bM u(x) \leq N\Lambda r^{-\sigma}\eta_{r,m}(x) \bM u(x).
\end{align*}
By using $\sum_m\eta_{r,m}^p\leq1$ and \eqref{eqmaximal},
\begin{align*}
\sum_{m\in\bZ^d} \|\eta_{r,m}L_{K_m}v_m\|_{L_{p,w}}^p &\leq N \Lambda^p r^{-\sigma p} \int_{\bR^d} \sum_{m\in\bZ^d} \eta_{r,m}(x)^p(\bM u(x))^pw(x)dx
\\
&\leq N\Lambda^pr^{-\sigma p}\|\bM u\|_{L_{p,w}}^p \leq N\Lambda^pr^{-\sigma p}\|u\|_{L_{p,w}}^p,
\end{align*}
which proves \eqref{eq11081405}.
The lemma is proved.
\end{proof}

\begin{remark}
    In Lemmas \ref{lem11131905} and \ref{lem11131906}, no regularity assumption on $K_m(x,y)$ in the $x$-variable is imposed, and only measurability is assumed.
\end{remark}

We next establish the boundedness of nonlocal operators defined by a family of localized kernels over a spatial grid. Recall that the truncated operator is introduced in \eqref{eq_trun}.

\begin{lemma} \label{lem11081808}
Let $\sigma\in(0,2)$, $p\in(1,\infty)$, $K_0>0$, $[w]_{A_p}\leq K_0$, $r\in(0,1]$, and $n_0\in\bN$. 
Suppose that $\Phi_{m,l}(x)$ $(m\in\bZ^d$ and $l=1,\ldots,n_0)$ are measurable functions such that $\Phi_{m,l}\geq0$, $\text{supp}(\Phi_{m,l})\subset B_r(x_{r,m})$, and
\begin{equation} \label{eq11090013}
\sum_{m\in\bZ^d} \sum_{l=1}^{n_0} \Phi_{m,l}(x)=1.
\end{equation}
For $m\in\bZ^d$ and $l=1,\ldots,n_0$, let $K_{m,l}(y)$ be measurable kernels satisfying Assumption \ref{as_bound} with the same constant $\Lambda$.
Then for $u\in C_c^\infty(\bR^d)$ and
\begin{equation} \label{eq11091343}
\overline{K}(x,y) := \sum_{m\in\bZ^d}\sum_{l=1}^{n_0} \Phi_{m,l}(x)K_{m,l}(y),
\end{equation}
we have
\begin{equation} \label{eq11082106}
\|L_{\overline{K},r}u\|_{L_{p,w}} \leq N\Lambda \|(-\Delta)^{\sigma/2}u\|_{L_{p,w}},
\end{equation}
where $N=N(d,p,K_0,\widetilde{\sigma},n_0)$ and $\widetilde{\sigma}$ is given by \eqref{eq_sigma}.
\end{lemma}

\begin{proof}
Let $\rho\in C_c^\infty(B_1)$ be a nonnegative radial function such that $\int_{\bR^d}\rho(x)dx=1$. Denote
\begin{equation*}
v_r(x) := \rho_r*u(x),
\end{equation*}
where $\rho_r(x):=r^{-d}\rho(x/r)$. Here, we emphasize that $r>0$ is fixed.
Then both $v_r$ and $u-v_r$ belong to $C_c^\infty(\bR^d)$.

Let $k=1$ if $\sigma\in(0,1)$ and $k=2$ if $\sigma\in[1,2)$.
We consider the Fourier multiplier symbols
\begin{equation*}
\frac{1-\hat{\rho}(r\xi)}{(r|\xi|)^\sigma}, \quad 1-\hat{\rho}(r\xi), \quad (r|\xi|)^{k-\sigma}\hat{\rho}(r\xi).
\end{equation*}
One can verify that they satisfy the Mikhlin condition \eqref{eq11082252}, with a constant $N(d)$ independent of $r>0$ and $\sigma\in(0,2)$. 
Since $\cF[v_r](\xi) = \hat{\rho}(r\xi)\hat{u}(\xi)$, by the weighted Mikhlin multiplier theorem and  \eqref{eq11050057}, we obtain
\begin{equation} \label{eq11082313}
r^{-\sigma}\|u-v_r\|_{L_{p,w}} + \|(-\Delta)^{\sigma/2}(u-v_r)\|_{L_{p,w}} \leq N \|(-\Delta)^{\sigma/2}u\|_{L_{p,w}}
\end{equation}
and
\begin{equation} \label{eq11082315}
r^{k-\sigma} \|D_x^kv_r\|_{L_{p,w}} \leq N \|(-\Delta)^{\sigma/2}u\|_{L_{p,w}},
\end{equation}
where $N=N(d,p,K_0)$. 

We first estimate $L_{\overline{K},r}v_r$. 
By a change of variables and \eqref{eq10291653} with $\gamma=k-\sigma$,
\begin{align*}
\int_{|y|<r} \frac{|D_x^kv_r(x+ty)|}{|y|^{d+\sigma-k}}dy &= t^{\sigma-k}\int_{|z|<tr} \frac{|D_x^kv_r(x+z)|}{|z|^{d+\sigma-k}} dz
\\
&\leq \frac{N}{k-\sigma}r^{k-\sigma} \bM(|D_x^kv_r|)(x).
\end{align*}
Since
\begin{equation*}
|\delta_y v_r(x)| \leq N |y|^k \int_0^1|D_x^kv_r(x+ty)|dt
\end{equation*}
for $|y|<r$, we get
\begin{equation*}
|L_{\overline{K},r}v_r(x)| \leq N\Lambda \frac{2-\sigma}{k-\sigma} r^{k-\sigma} \bM(|D_x^kv_r|)(x).
\end{equation*}
Therefore, \eqref{eqmaximal} and \eqref{eq11082315} imply
\begin{equation} \label{eq11090022}
\|L_{\overline{K},r}v_r\|_{L_{p,w}} \leq N \Lambda \|(-\Delta)^{\sigma/2}u\|_{L_{p,w}}.
\end{equation}

It remains to estimate $L_{\overline{K},r}(u-v_r)$.
If $\Phi_{m,l}(x)\ne0$, then $x\in B_r(x_{r,m})$, so $x,x+y\in B_{2r}(x_{r,m})$ whenever $|y|<r$. Thus, using $\widetilde{\eta}_{r,m}=1$ on $B_{4r}(x_{r,m})$, we deduce that for $|y|<r$,
\begin{equation*}
\delta_y(u-v_r)(x)=\delta_y(\widetilde{\eta}_{r,m}(u-v_r))(x).
\end{equation*}
 Hence,
\begin{equation*}
L_{\overline{K},r}(u-v_r)(x)=\sum_{m\in\bZ^d}\sum_{l=1}^{n_0} \Phi_{m,l}(x) L_{K_{m,l},r}(\widetilde{\eta}_{r,m}(u-v_r))(x).
\end{equation*}
Since $K_{m,l}$ satisfies Assumption \ref{as_bound}, by \eqref{eq10311722} and \eqref{eq11090013},
\begin{align*}
\|L_{\overline{K},r}(u-v_r)\|_{L_{p,w}}^p &\leq \sum_{m\in\bZ^d}\sum_{l=1}^{n_0} \|L_{K_{m,l},r}(\widetilde{\eta}_{r,m}(u-v_r))\|_{L_{p,w}}^p
\\
&\leq N n_0 \Lambda^p \sum_{m\in\bZ^d}\|(-\Delta)^{\sigma/2}(\widetilde{\eta}_{r,m}(u-v_r))\|_{L_{p,w}}^p.
\end{align*}
By \eqref{eq11081703} with $\zeta_m=\widetilde{\eta}_{r,m}$ and $\kappa=8$, and \eqref{eq11082313},
\begin{align} \label{eq11090023}
\|L_{\overline{K},r}(u-v_r)\|_{L_{p,w}}&\leq N\Lambda\left(\sum_m\|(-\Delta)^{\sigma/2}(\widetilde{\eta}_{r,m}(u-v_r))\|_{L_{p,w}}^p\right)^{1/p} \nonumber
\\
&\leq N\Lambda\left(\|(-\Delta)^{\sigma/2}(u-v_r)\|_{L_{p,w}}+r^{-\sigma}\|u-v_r\|_{L_{p,w}}\right) \nonumber
\\
&\leq N\Lambda\|(-\Delta)^{\sigma/2}u\|_{L_{p,w}}.
\end{align}

Due to $u=v_r+(u-v_r)$, \eqref{eq11082106} follows from \eqref{eq11090022} and \eqref{eq11090023}.
The lemma is proved.
\end{proof}

For the remainder of this section, we replace $\omega(s)$ by $\min\{\omega(s),2\Lambda\}$ for $0<s\leq1$. Then under Assumption \ref{as_bound}, we still have \eqref{eq_dini}. 
We recall that
\begin{equation*}
N_\omega := \int_0^1\omega(r)\frac{dr}{r}.
\end{equation*}

For $x_0\in\bR^d$, we write $K_{x_0}(y):=K(x_0,y)$.
Denote $\chi_{r,x_0}(x):=\widetilde{\eta}((x-x_0)/r)$, where $\widetilde{\eta}$ is introduced in \eqref{eq11090037}. Thus, $\chi_{r,x_0}=1$ on $B_{2r}(x_0)$ and is supported in $B_{4r}(x_0)$.

We next present a local estimate for the difference between $L_K$ and $L_{K_{x_0}}$.

\begin{lemma}
Let $\sigma\in(0,2)$, $p\in(1,\infty)$, $K_0>0$, and $[w]_{A_p}\leq K_0$. Suppose that $K$ satisfies Assumptions \ref{as_bound} and \ref{as_dini}.
Then for $0<r\leq1/16$, $x_0\in\bR^d$, and $u\in C_c^\infty(\bR^d)$,
\begin{equation} \label{eq11090054}
\|\chi_{r,x_0}(L_K-L_{K_{x_0}})u\|_{L_{p,w}} \leq N N_\omega(16r) \|(-\Delta)^{\sigma/2}u\|_{L_{p,w}} + N\Lambda r^{-\sigma} \|u\|_{L_{p,w}},
\end{equation}
where $N=N(d,p,K_0,\widetilde{\sigma})$ is independent of $x_0,r,\Lambda$, and $\omega$.
\end{lemma}

\begin{proof}
\textbf{Step 1.} First, we construct a nonnegative kernel and show that its increments are controlled by the original modulus $\omega$. Define
\begin{equation} \label{eq11091800}
\widetilde{K}(x,y) := (2-\sigma)\omega(8r) + \chi_{2r,x_0}(x)(K(x,y)-K_{x_0}(y)).
\end{equation}
Since $\text{supp}(\chi_{2r,x_0})\subset B_{8r}(x_0)$ and $8r\leq1$, by \eqref{eq_dini},
\begin{equation*}
|\chi_{2r,x_0}(x)(K(x,y)-K_{x_0}(y))|\leq(2-\sigma)\omega(8r),
\end{equation*}
which yields
\begin{equation} \label{eq11091355}
0\leq\widetilde{K}(x,y) \leq 2(2-\sigma) \omega(8r).
\end{equation}
Moreover, when $\sigma=1$, one can see that \eqref{eq_cancel} holds for $\widetilde{K}$.

Let $0<s\leq8r$ and $0<|h|\leq s$. We use the identity
\begin{align} \label{eq11140107}
\widetilde{K}(x+h,y)-\widetilde{K}(x,y) &= \chi_{2r,x_0}(x+h)(K(x+h,y)-K(x,y)) \nonumber
\\
&\quad+ (\chi_{2r,x_0}(x+h)-\chi_{2r,x_0}(x)) (K(x,y)-K_{x_0}(y)).
\end{align}
Due to $s\leq8r\leq1$ and \eqref{eq_dini},
\begin{equation} \label{eq11091332}
|\chi_{2r,x_0}(x+h)(K(x+h,y)-K(x,y))|\leq (2-\sigma)\omega(|h|) \leq (2-\sigma)\omega(s).
\end{equation}
For the second term on the right-hand side of \eqref{eq11140107}, if $\chi_{2r,x_0}(x+h)-\chi_{2r,x_0}(x)\neq0$, then either $x$ or $x+h$ belongs to $B_{8r}(x_0)$, and hence
\begin{equation*}
|x-x_0|\leq8r+|h|\leq8r+s\leq16r.
\end{equation*}
Divide the segment joining $x_0$ and $x$ into $\left\lceil\frac{16r}{s}\right\rceil$ equal subsegments so that each has length at most $s$. Then by applying \eqref{eq_dini} along these subsegments, one can deduce that
\begin{equation*}
|K(x,y)-K_{x_0}(y)| \leq \left\lceil\frac{16r}{s}\right\rceil(2-\sigma) \omega(s) \leq \frac{32r}{s} (2-\sigma) \omega(s).
\end{equation*}
Using $|D_x\chi_{2r,x_0}|\leq\frac{N}{r}$,
\begin{align*}
|\chi_{2r,x_0}(x+h)-\chi_{2r,x_0}(x)||K(x,y)-K_{x_0}(y)| &\leq N \frac{s}{r}\frac{32r}{s}(2-\sigma)\omega(s)
\\
&\leq N(2-\sigma)\omega(s).
\end{align*}
Therefore, this and \eqref{eq11091332} yield that for $0<|h|\leq s\leq8r$,
\begin{equation} \label{eq11091333}
|\widetilde{K}(x+h,y)-\widetilde{K}(x,y)| \leq N(d) (2-\sigma)\omega(s).
\end{equation}

\textbf{Step 2.}
In this step, we establish an approximation of $\widetilde{K}$.
Let $r_j:=2^{-j}$ for $j\in\bN_0$, and choose $j_0\in\bN$ so that $r\leq2^{-j_0}<2r$.
For $k\geq j_0$, define
\begin{equation} \label{eq11091426}
\widetilde{K}_k(x,y) := \sum_{m\in\bZ^d} \eta_{r_k,m}(x) \widetilde{K}(x_{r_k,m},y),
\end{equation}
where $\eta_{r_k,m}$ is defined in \eqref{eqeta}.
By \eqref{eq11091355}, $0 \leq \widetilde{K}_k \leq 2(2-\sigma)\omega(8r)$.
In particular, $\widetilde{K}_k$ satisfies Assumption \ref{as_bound} with $2\omega(8r)$ in place of $\Lambda$.
Since $\eta_{r_k,m}(x)\neq0$ only if $|x-x_{r_k,m}|<r_k<2r$, by \eqref{eq11091424} and \eqref{eq11091333},
\begin{equation} \label{eq11092149}
|\widetilde{K}(x,y)-\widetilde{K}_k(x,y)| \leq N(2-\sigma)\omega(r_k).
\end{equation}

We also show that the difference between the approximating kernels $\widetilde{K}_j(x,y)-\widetilde{K}_{j-1}(x,y)$ can be represented by kernels of the form \eqref{eq11091343}.
For $j>j_0$ and $q\in\bZ^d$, if $\eta_{r_j,m}(x)\eta_{r_{j-1},q}(x)\neq0$, then
\begin{equation*}
|x_{r_j,m}-x_{r_{j-1},q}| < r_j+r_{j-1} < 4r_j,
\end{equation*}
which yields $|m-2q|<4\sqrt{d}$.
For each $m\in\bZ^d$, let $q_{m,1}, \ldots, q_{m,n_m}$ be all the elements of $\bZ^d$ satisfying $|m-2q|<4\sqrt{d}$.
Here, we see that $n_m\leq(4\sqrt{d}+1)^d$.

Let $n_0=n_0(d)$ be any integer satisfying $n_0\geq (4\sqrt{d}+1)^d$.
Denote
\begin{equation*}
\Phi^{(j)}_{m,l}(x) := \eta_{r_j,m}(x) \eta_{r_{j-1},q_{m,l}}(x), \quad 1\leq l\leq n_m,
\end{equation*}
and let $\Phi^{(j)}_{m,l}=0$ for $n_m<l\leq n_0$. 
Using \eqref{eq11091424},
\begin{equation} \label{eq11092212}
\Phi^{(j)}_{m,l}\geq0, \quad \text{supp}(\Phi^{(j)}_{m,l}) \subset B_{r_j}(x_{r_j,m}), \quad \sum_{m\in\bZ^d} \sum_{l=1}^{n_0} \Phi^{(j)}_{m,l}(x) = 1.
\end{equation}
Consequently, by \eqref{eq11091426},
\begin{align} \label{eq11091434}
&\widetilde{K}_j(x,y)-\widetilde{K}_{j-1}(x,y) \nonumber
\\
&= \sum_{m\in\bZ^d} \sum_{q\in\bZ^d} \eta_{r_j,m}(x)\eta_{r_{j-1},q}(x) (\widetilde{K}(x_{r_j,m},y)-\widetilde{K}(x_{r_{j-1},q},y)) \nonumber
\\
&= \sum_{m\in\bZ^d} \sum_{l=1}^{n_m} \Phi^{(j)}_{m,l}(x) (\widetilde{K}(x_{r_j,m},y)-\widetilde{K}(x_{r_{j-1},q_{m,l}},y)). 
\end{align}

\textbf{Step 3.} 
In this step, our goal is to prove
\begin{align} \label{eq11091757}
\|\chi_{r,x_0}(L_K-L_{K_{x_0}})u\|_{L_{p,w}} &\leq \|L_{\widetilde{K},2^{-j_0}}u\|_{L_{p,w}} + N N_\omega(16r) \|(-\Delta)^{\sigma/2}u\|_{L_{p,w}} \nonumber
\\
&\quad+ N\Lambda r^{-\sigma}\|u\|_{L_{p,w}}.
\end{align}
Here, the operator $L_{\widetilde{K},2^{-j_0}}$ is defined as in \eqref{eq_trun}.

Since $\chi_{2r,x_0}=1$ on $\text{supp}\,\chi_{r,x_0}$, due to \eqref{eq11091800},
\begin{align} \label{eq11091848}
\chi_{r,x_0}(x)(L_K-L_{K_{x_0}})u(x) &= \chi_{r,x_0}(x)L_{\widetilde{K},2^{-j_0}}u(x)-\chi_{r,x_0}(x)L_{(2-\sigma)\omega(8r),2^{-j_0}}u(x) \nonumber
\\
&\quad+ \chi_{r,x_0}(x)\int_{|y|\geq2^{-j_0}}\delta_yu(x)(K(x,y)-K_{x_0}(y))\frac{dy}{|y|^{d+\sigma}}.
\end{align}
Thus, to obtain \eqref{eq11091757}, we need to estimate the second and third terms on the right-hand side.

Applying \eqref{eq10311722} to the constant kernel $(2-\sigma)\omega(8r)$, we obtain
\begin{align} \label{eq11091849}
\|\chi_{r,x_0}L_{(2-\sigma)\omega(8r),2^{-j_0}}u\|_{L_{p,w}} &\leq N\omega(8r) \|(-\Delta)^{\sigma/2}u\|_{L_{p,w}} \nonumber
\\
&\leq N N_\omega(16r) \|(-\Delta)^{\sigma/2}u\|_{L_{p,w}}.
\end{align}
Here, for the last inequality, we used $\omega(8r)\leq(\log{2})^{-1}N_\omega(16r)$, which follows from \eqref{eq10301427}.

For the third term, by \eqref{eq_dini}, for $x\in\text{supp}(\chi_{r,x_0})$,
\begin{equation*}
|K(x,y)-K_{x_0}(y)|\leq(2-\sigma)\omega(4r).
\end{equation*}
Thus, by \eqref{eq_cancel} and \eqref{eq10291650},
\begin{align*}
&\left| \chi_{r,x_0}(x)\int_{|y|\geq2^{-j_0}} \delta_yu(x)(K(x,y)-K_{x_0}(y)) \frac{dy}{|y|^{d+\sigma}} \right|
\\
&\leq (2-\sigma)\omega(4r) \int_{|y|\geq2^{-j_0}} \left(\frac{|u(x+y)|+|u(x)|}{|y|^{d+\sigma}} + 1_{\{\sigma>1\}} \frac{|\nabla u(x)|}{|y|^{d+\sigma-1}} \right) dy
\\
&\leq N \frac{2-\sigma}{\sigma} \omega(4r)2^{j_0\sigma}(\bM u(x)+|u(x)|) + N 1_{\{\sigma>1\}} \frac{2-\sigma}{\sigma-1} \omega(4r) 2^{j_0(\sigma-1)}|D_xu(x)|
\\
&\leq N\Lambda r^{-\sigma}(\bM u(x)+|u(x)|) + N 1_{\{\sigma>1\}} \omega(4r) r^{1-\sigma}|D_xu(x)|.
\end{align*}
For the last inequality, we used $2^{-j_0}\geq r$ and $\omega(4r) \leq 2 \Lambda$.
When $\sigma>1$, by applying \eqref{eq10301427}, and \eqref{eq_grad} with $\varepsilon=r$,
\begin{align*}
\omega(4r) r^{1-\sigma} \|D_xu\|_{L_{p,w}} &\leq N\omega(4r)\|(-\Delta)^{\sigma/2}u\|_{L_{p,w}} + N\omega(4r) r^{-\sigma} \|u\|_{L_{p,w}}
\\
&\leq N N_\omega(8r) \|(-\Delta)^{\sigma/2}u\|_{L_{p,w}} + N\Lambda r^{-\sigma}\|u\|_{L_{p,w}}.
\end{align*}
Hence, by \eqref{eqmaximal},
\begin{align*}
&\left\|\chi_{r,x_0}\int_{|y|\geq2^{-j_0}}\delta_yu(K-K_{x_0})\frac{dy}{|y|^{d+\sigma}}\right\|_{L_{p,w}}
\\
&\leq NN_\omega(8r)\|(-\Delta)^{\sigma/2}u\|_{L_{p,w}} + N\Lambda r^{-\sigma}\|u\|_{L_{p,w}}.
\end{align*}
This together with \eqref{eq11091848} and \eqref{eq11091849} yields \eqref{eq11091757}.

\textbf{Step 4.} 
To conclude \eqref{eq11090054}, it remains to bound $\|L_{\widetilde{K},2^{-j_0}}u\|_{L_{p,w}}$ in \eqref{eq11091757}.
In this step, we decompose $L_{\widetilde{K},2^{-j_0}}u$.
For $n\geq j_0$, we write
\begin{align} \label{eq11092134}
L_{\widetilde{K},2^{-j_0}}u(x) &= L_{\widetilde{K},r_{n+1}}u(x) + \sum_{k=j_0}^{n} \int_{\{r_{k+1}\leq|y|<r_k\}} \delta_yu(x) \widetilde{K}(x,y) \frac{dy}{|y|^{d+\sigma}} \nonumber
\\
&= L_{\widetilde{K},r_{n+1}}u(x) + \sum_{k=j_0}^{n} \int_{\{r_{k+1}\leq|y|<r_k\}} \delta_yu(x) (\widetilde{K}-\widetilde{K}_k)(x,y) \frac{dy}{|y|^{d+\sigma}} \nonumber
\\
&\quad+ \sum_{k=j_0}^{n} \int_{\{r_{k+1}\leq|y|<r_k\}} \delta_yu(x) \widetilde{K}_k(x,y) \frac{dy}{|y|^{d+\sigma}}.
\end{align}
Here, we used the identity $\widetilde{K}=(\widetilde{K}-\widetilde{K}_k)+\widetilde{K}_k$.

For the last term above, by using $\widetilde{K}_k=\widetilde{K}_{j_0}+\sum_{j=j_0+1}^{k}(\widetilde{K}_j-\widetilde{K}_{j-1})$, we obtain
\begin{align} \label{eq11092135}
&\sum_{k=j_0}^{n} \int_{\{r_{k+1}\leq|y|<r_k\}} \delta_yu(x) \widetilde{K}_k(x,y) \frac{dy}{|y|^{d+\sigma}} \nonumber
\\
&= \int_{\{r_{n+1}\leq|y|<2^{-j_0}\}} \delta_yu(x) \widetilde{K}_{j_0}(x,y) \frac{dy}{|y|^{d+\sigma}} \nonumber
\\
&\quad+ \sum_{k=j_0+1}^{n} \sum_{j=j_0+1}^{k} \int_{\{r_{k+1}\leq|y|<r_k\}} \delta_yu(x) (\widetilde{K}_j-\widetilde{K}_{j-1})(x,y) \frac{dy}{|y|^{d+\sigma}} \nonumber
\\
&= \int_{\{r_{n+1}\leq|y|<2^{-j_0}\}} \delta_yu(x) \widetilde{K}_{j_0}(x,y) \frac{dy}{|y|^{d+\sigma}} \nonumber
\\
&\quad+ \sum_{j=j_0+1}^{n} \int_{\{r_{n+1}\leq|y|<r_j\}} \delta_yu(x)(\widetilde{K}_j-\widetilde{K}_{j-1})(x,y) \frac{dy}{|y|^{d+\sigma}} \nonumber
\\
&=L_{\widetilde{K}_{j_0},2^{-j_0}}u(x) -L_{\widetilde{K}_{j_0},r_{n+1}}u(x) \nonumber
\\
&\quad+ \sum_{j=j_0+1}^{n}(L_{\widetilde{K}_j,r_j}-L_{\widetilde{K}_{j-1},r_j})u(x) - \sum_{j=j_0+1}^{n}(L_{\widetilde{K}_j,r_{n+1}} - L_{\widetilde{K}_{j-1},r_{n+1}})u(x).
\end{align}
Since
\begin{equation*}
L_{\widetilde{K}_{j_0},r_{n+1}}u + \sum_{j=j_0+1}^{n}(L_{\widetilde{K}_j,r_{n+1}} - L_{\widetilde{K}_{j-1},r_{n+1}})u = L_{\widetilde{K}_n,r_{n+1}}u,
\end{equation*}
combining \eqref{eq11092134} and \eqref{eq11092135} yields
\begin{align} \label{eq11092136}
L_{\widetilde{K},2^{-j_0}}u(x) &= \sum_{k=j_0}^{n} \int_{\{r_{k+1}\leq|y|<r_k\}} \delta_yu(x) (\widetilde{K}-\widetilde{K}_k)(x,y) \frac{dy}{|y|^{d+\sigma}} \nonumber
\\
&\quad+ L_{\widetilde{K}_{j_0},2^{-j_0}}u(x) +\sum_{j=j_0+1}^{n}(L_{\widetilde{K}_j,r_j}-L_{\widetilde{K}_{j-1},r_j})u(x) \nonumber
\\
&\quad+(L_{\widetilde{K},r_{n+1}}-L_{\widetilde{K}_n,r_{n+1}})u(x).
\end{align}

\textbf{Step 5.}
In this step, we estimate the first three terms on the right-hand side of \eqref{eq11092136}.
Let us consider the first term. By \eqref{eq_shell} and \eqref{eq11092149},
\begin{align*}
&\left\|\int_{\{r_{k+1}\leq|y|<r_k\}} \delta_yu(\widetilde{K}-\widetilde{K}_k) \frac{dy}{|y|^{d+\sigma}}\right\|_{L_{p,w}}
\\
&\leq N(2-\sigma)\omega(r_k) \left\|\int_{\{r_{k+1}\leq|y|<r_k\}}\frac{|\delta_yu|}{|y|^{d+\sigma}}dy\right\|_{L_{p,w}} \leq N \omega(r_k) \|(-\Delta)^{\sigma/2}u\|_{L_{p,w}}.
\end{align*}
Since by \eqref{eq10301427},
\begin{align*}
\sum_{k\geq j_0}\omega(4r_k) &\leq\frac{1}{\log{2}}N_\omega(2^{3-j_0})\leq\frac{1}{\log{2}}N_\omega(16r),
\end{align*}
we have
\begin{align} \label{eq11092237}
&\sum_{k=j_0}^{n} \left\|\int_{\{r_{k+1}\leq|y|<r_k\}} \delta_yu(\widetilde{K}-\widetilde{K}_k) \frac{dy}{|y|^{d+\sigma}}\right\|_{L_{p,w}} \leq N N_\omega(16r) \|(-\Delta)^{\sigma/2}u\|_{L_{p,w}}.
\end{align}

For the second term, we apply Lemma \ref{lem11081808} with $2^{-j_0}$ in place of $r$, $n_0=1$, $\Phi_{m,1}=\eta_{2^{-j_0},m}$, the kernels $\widetilde{K}(x_{2^{-j_0},m},y)$, and $2\omega(8r)$ in place of $\Lambda$. By \eqref{eq10301427}, \eqref{eq11091355}, and \eqref{eq11091426},
\begin{equation} \label{eq11092238}
\|L_{\widetilde{K}_{j_0},2^{-j_0}}u\|_{L_{p,w}}\leq N\omega(8r)\|(-\Delta)^{\sigma/2}u\|_{L_{p,w}} \leq NN_\omega(16r)\|(-\Delta)^{\sigma/2}u\|_{L_{p,w}}.
\end{equation}

For the third term, we use \eqref{eq11091434} and add a constant to make each frozen kernel nonnegative. For $j>j_0$, $m\in\bZ^d$, and $1\leq l\leq n_m$, define
\begin{equation*}
K^{(j)}_{m,l}(y) := \widetilde{K}(x_{r_j,m},y) - \widetilde{K}(x_{r_{j-1},q_{m,l}},y) + c_j,
\end{equation*}
where
\begin{equation*}
c_j:=C_0(2-\sigma)\omega(4r_j),
\end{equation*}
and let $K^{(j)}_{m,l}:=0$ for $n_m<l\leq n_0$.
Here, due to \eqref{eq11091333}, one can choose $C_0=C_0(d)>0$ sufficiently large so that
\begin{equation*}
0\leq K^{(j)}_{m,l}(y)\leq N(2-\sigma)\omega(4r_j).
\end{equation*}

Let
\begin{equation*}
\overline{K}_j(x,y) := \sum_{m\in\bZ^d} \sum_{l=1}^{n_0} \Phi^{(j)}_{m,l}(x)K^{(j)}_{m,l}(y).
\end{equation*}
Due to \eqref{eq11092212} and \eqref{eq11091434},
\begin{align*}
\overline{K}_j&=\widetilde{K}_j-\widetilde{K}_{j-1}+c_j.
\end{align*}
By Lemma \ref{lem11081808} with $N\omega(4r_j)$ in place of $\Lambda$, and applying \eqref{eq10311722} to the constant kernel $c_j$,
\begin{align} \label{eq11092239}
\|(L_{\widetilde{K}_j,r_j}-L_{\widetilde{K}_{j-1},r_j})u\|_{L_{p,w}} &\leq \|L_{\overline{K}_j,r_j}u\|_{L_{p,w}} + \|L_{c_j,r_j}u\|_{L_{p,w}} \nonumber
\\
&\leq N\omega(4r_j) \|(-\Delta)^{\sigma/2}u\|_{L_{p,w}}.
\end{align}

\textbf{Step 6.} 
Now we prove \eqref{eq11090054}.
By \eqref{eq11092136}, \eqref{eq11092237}, \eqref{eq11092238}, and \eqref{eq11092239}, for any integer $n\geq j_0$,
\begin{align*}
\|(L_{\widetilde{K},2^{-j_0}} - L_{\widetilde{K},r_{n+1}} + L_{\widetilde{K}_n,r_{n+1}})u\|_{L_{p,w}} \leq N N_\omega(16r) \|(-\Delta)^{\sigma/2}u\|_{L_{p,w}}.
\end{align*}
Due to \eqref{eq11092149},
\begin{align*}
|(L_{\widetilde{K},r_{n+1}}-L_{\widetilde{K}_n,r_{n+1}})u(x)| &\leq N(2-\sigma)\omega(r_n) \int_{|y|<r_{n+1}} \frac{|\delta_yu(x)|}{|y|^{d+\sigma}}dy
\\
&\leq N(2-\sigma)\omega(r_n) \int_{|y|<1} \frac{|\delta_yu(x)|}{|y|^{d+\sigma}} dy\to0
\end{align*}
as $n\to\infty$.
By Fatou's lemma, we conclude that
\begin{equation*}
\|L_{\widetilde{K},2^{-j_0}}u\|_{L_{p,w}} \leq N N_\omega(16r) \|(-\Delta)^{\sigma/2}u\|_{L_{p,w}}.
\end{equation*}
This together with \eqref{eq11091757} yields \eqref{eq11090054}.
The lemma is proved.
\end{proof}

\section{Proof of the main results} \label{sec_proof}

\subsection{Proof of the boundedness of operators}

\begin{proof}[Proof of Theorem \ref{thm_bound}]

By Remark \ref{rem_dense}, we only need to prove \eqref{est_bound} when $u\in C_c^\infty(\bR^d)$.

Let $r=\frac{1}{16}$ and $K_m(y):=K(x_{r,m},y)$ for $m\in\bZ^d$. 
Note that
\begin{align*}
\eta_{r,m}L_Ku = \eta_{r,m}L_{K_m}(\widetilde{\eta}_{r,m}u) + \eta_{r,m}(L_K-L_{K_m})(\widetilde{\eta}_{r,m}u) + \eta_{r,m}L_K((1-\widetilde{\eta}_{r,m})u),
\end{align*}
where $\eta_{r,m}$ and $\widetilde{\eta}_{r,m}$ are taken from \eqref{eqeta}.
By \eqref{eq10311722},
\begin{equation*}
\|\eta_{r,m}L_{K_m}(\widetilde{\eta}_{r,m}u)\|_{L_{p,w}} \leq N\Lambda\|(-\Delta)^{\sigma/2}(\widetilde{\eta}_{r,m}u)\|_{L_{p,w}}.
\end{equation*}
Let $\chi_{r,x_{r,m}}(x):=\widetilde{\eta}((x-x_{r,m})/r)$, where $\widetilde{\eta}$ is introduced in \eqref{eq11090037}.
Since $\chi_{r,x_{r,m}}=1$ on $\text{supp}(\eta_{r,m})$, applying \eqref{eq11090054} with $x_0=x_{r,m}$ and $\widetilde{\eta}_{r,m}u$ in place of $u$,
\begin{align} \label{eq11100122}
&\|\eta_{r,m}(L_K-L_{K_m})(\widetilde{\eta}_{r,m}u)\|_{L_{p,w}} \nonumber
\\
&\leq \|\chi_{r,x_{r,m}}(L_K-L_{K_m})(\widetilde{\eta}_{r,m}u)\|_{L_{p,w}} \nonumber
\\
&\leq NN_\omega(16r) \|(-\Delta)^{\sigma/2}(\widetilde{\eta}_{r,m}u)\|_{L_{p,w}} + N\Lambda r^{-\sigma} \|\widetilde{\eta}_{r,m}u\|_{L_{p,w}}.
\end{align}
By \eqref{eq11081405},
\begin{equation*}
\left(\sum_{m\in\bZ^d} \|\eta_{r,m}L_K((1-\widetilde{\eta}_{r,m})u)\|_{L_{p,w}}^p\right)^{1/p} \leq N \Lambda r^{-\sigma} \|u\|_{L_{p,w}}.
\end{equation*}
Thus, by \eqref{eq11061351}, \eqref{eq11061352}, and \eqref{eq11081703} with $\zeta_m=\widetilde{\eta}_{r,m}$ and $\kappa=8$,
\begin{align*}
\|L_Ku\|_{L_{p,w}} &\leq N\left(\sum_{m\in\bZ^d}\|\eta_{r,m}L_Ku\|_{L_{p,w}}^p\right)^{1/p}
\\
&\leq N(\Lambda+N_\omega(16r)) \left(\sum_{m\in\bZ^d}\|(-\Delta)^{\sigma/2}(\widetilde{\eta}_{r,m}u)\|_{L_{p,w}}^p\right)^{1/p}
\\
&\quad+ N\Lambda r^{-\sigma}\left(\sum_{m\in\bZ^d} \|\widetilde{\eta}_{r,m}u\|_{L_{p,w}}^p\right)^{1/p} + N\Lambda r^{-\sigma}\|u\|_{L_{p,w}}
\\
&\leq N(\Lambda+N_\omega(16r)) \left(\|(-\Delta)^{\sigma/2}u\|_{L_{p,w}} + r^{-\sigma}\|u\|_{L_{p,w}}\right).
\end{align*}
Since $r=\frac{1}{16}$ is fixed, by \eqref{eq10310040} and $N_\omega(16r)=N_{\omega}$,
\begin{equation*}
\|L_Ku\|_{L_{p,w}} \leq N(\Lambda+N_\omega)\|u\|_{H_{p,w}^\sigma},
\end{equation*}
where $N=N(d,p,K_0,\widetilde{\sigma})$.
The theorem is proved.
\end{proof}

\subsection{Proof of the a priori estimate and solvability}

\begin{proof}[Proof of Theorem \ref{thm_dini}]

Recall that the solvability for kernels independent of $x$ is already obtained in Lemma \ref{lem_simple}. Thus, due to the method of continuity, it suffices to prove \eqref{est_main}. Moreover, by Theorem \ref{thm_bound}, we may assume that $u\in C_c^\infty(\bR^d)$.

Let $\lambda>0$ and $f:=L_Ku-\lambda u$.
Let $r\in(0,\frac{1}{16}]$, which will be determined later.
Note that for each $m\in\bZ^d$, $K_m(y):=K(x_{r,m},y)$ satisfies Assumptions \ref{as_bound} and \ref{as_nonde} with the same constants $\Lambda$ and $\nu$. Thus, by \eqref{eq10311738},
\begin{equation*}
\|(-\Delta)^{\sigma/2}(\eta_{r,m}u)\|_{L_{p,w}}+\lambda\|\eta_{r,m}u\|_{L_{p,w}} \leq N\|(\lambda-L_{K_m})(\eta_{r,m}u)\|_{L_{p,w}}.
\end{equation*}
By \eqref{eq11081702} and \eqref{eq11061351}, we obtain
\begin{align} \label{eq11100125}
&\|(-\Delta)^{\sigma/2}u\|_{L_{p,w}} + \lambda\|u\|_{L_{p,w}} \nonumber
\\
&\leq N \left(\sum_{m\in\bZ^d} \|(-\Delta)^{\sigma/2}(\eta_{r,m}u)\|_{L_{p,w}}^p\right)^{1/p} \nonumber
\\
&\quad+ N \lambda \left(\sum_{m\in\bZ^d}\|\eta_{r,m}u\|_{L_{p,w}}^p\right)^{1/p} + Nr^{-\sigma}\|u\|_{L_{p,w}} \nonumber
\\
&\leq N \left(\sum_{m\in\bZ^d} \|(\lambda-L_{K_m})(\eta_{r,m}u)\|_{L_{p,w}}^p\right)^{1/p} + N r^{-\sigma} \|u\|_{L_{p,w}}.
\end{align}

Now we consider
\begin{align*} 
(\lambda-L_{K_m})(\eta_{r,m}u) &= -\eta_{r,m}f+\eta_{r,m}(L_K-L_{K_m})(\widetilde{\eta}_{r,m}u)
\\
&\quad+\eta_{r,m}L_K((1-\widetilde{\eta}_{r,m})u)-\eta_{r,m}L_{K_m}((1-\widetilde{\eta}_{r,m})u)
\\
&\quad-[L_{K_m},\eta_{r,m}]u.
\end{align*}
Here, $[\cdot,\cdot]$ is defined in \eqref{eq11081725}.
For the first term on the right-hand side, due to \eqref{eq11061351},
\begin{equation*}
\left(\sum_{m\in\bZ^d} \|\eta_{r,m}f\|_{L_{p,w}}^p\right)^{1/p}\leq \|f\|_{L_{p,w}}.
\end{equation*}
The second term is the same as that in \eqref{eq11100122} and hence satisfies the same bound.
For the two terms involving $(1-\widetilde{\eta}_{r,m})u$, we apply \eqref{eq11081405} separately to the family consisting of the same kernel $K$ and to the family of frozen kernels $K_m$. This gives
\begin{align*}
&\left(\sum_{m\in\bZ^d}\|\eta_{r,m}L_K((1-\widetilde{\eta}_{r,m})u)\|_{L_{p,w}}^p\right)^{1/p}
\\
&\quad+\left(\sum_{m\in\bZ^d}\|\eta_{r,m}L_{K_m}((1-\widetilde{\eta}_{r,m})u)\|_{L_{p,w}}^p\right)^{1/p} \leq N\Lambda r^{-\sigma}\|u\|_{L_{p,w}}.
\end{align*}
For the last term, by \eqref{eq11071436} with $\zeta_m=\eta_{r,m}$ and $\kappa=1$, for $\varepsilon\in(0,1)$,
\begin{equation*}
\left(\sum_{m\in\bZ^d} \|[L_{K_m},\eta_{r,m}]u\|_{L_{p,w}}^p\right)^{1/p} \leq N \Lambda \varepsilon\|(-\Delta)^{\sigma/2}u\|_{L_{p,w}} + N_\varepsilon \Lambda r^{-\sigma}\|u\|_{L_{p,w}}.
\end{equation*}

Combining these estimates with \eqref{eq11100125},
\begin{align} \label{eq11100126}
\|(-\Delta)^{\sigma/2}u\|_{L_{p,w}} + \lambda\|u\|_{L_{p,w}} &\leq N \|f\|_{L_{p,w}} + N(N_\omega(16r)+\Lambda\varepsilon) \|(-\Delta)^{\sigma/2}u\|_{L_{p,w}} \nonumber
\\
&\quad+ N_\varepsilon(1+\Lambda+N_\omega)r^{-\sigma} \|u\|_{L_{p,w}}.
\end{align}

Since $N_\omega(s)\to0$ as $s\to0$, we can choose $\varepsilon>0$ and $r\in(0,\frac{1}{16}]$ sufficiently small so that
\begin{equation*}
N(N_\omega(16r)+\Lambda\varepsilon) \leq \frac{1}{2}.
\end{equation*}
Then we denote
\begin{equation*}
\lambda_0:=2N_\varepsilon(1+\Lambda+N_\omega)r^{-2}.
\end{equation*}
Thus, for $\lambda\geq\lambda_0$, by \eqref{eq11100126} and $r^{-\sigma}\leq r^{-2}$,
\begin{align*}
\|(-\Delta)^{\sigma/2}u\|_{L_{p,w}} + \lambda\|u\|_{L_{p,w}} &\leq N \|f\|_{L_{p,w}} + \frac{1}{2}\|(-\Delta)^{\sigma/2}u\|_{L_{p,w}} + \frac{\lambda}{2}\|u\|_{L_{p,w}}.
\end{align*}
Now \eqref{est_main} can be obtained by absorbing the last two terms of the right-hand side.
The theorem is proved.
\end{proof}

\section{Global counterexamples with uniformly continuous kernels} \label{sec_counter}

We prove Theorem \ref{thm_counter} by adapting the construction in \cite{L26}.

\subsection{Construction of the kernel} \label{sec_cons}

Let $\phi\in C_c^\infty(\bR^d)$ be a nonnegative radial function such that
\begin{equation} \label{eq11101439}
\text{supp}(\phi) \subset \left\{\xi\in\bR^d:\frac{1}{2}<|\xi|<2\right\}, \quad \phi(e_1)=1,
\end{equation}
where $e_1=(1,0,\ldots,0)$.
 Define
\begin{equation*}
\hat{q}(\xi):=|\xi|^\sigma\phi(\xi), \quad q_j(y):=2^{j(d+\sigma)}q(2^jy), \quad j\geq1.
\end{equation*}
Then $q$ is a real valued radial function in $\cS(\bR^d)$. Moreover,
\begin{equation*}
\int_{\bR^d}q(y)dy=0, \quad \hat{q}_j(\xi) = |\xi|^\sigma\phi(2^{-j}\xi).
\end{equation*}
Since $q\in \cS(\bR^d)$, for any $M\in\bN_0$,
\begin{align*}
|q(y)| &= \frac{1}{(2\pi)^d(1+|y|^2)^M}\left|\int_{\bR^d}e^{iy\cdot\xi}(1-\Delta_\xi)^M\left(|\xi|^\sigma\phi(\xi)\right)d\xi\right|
\\
&\leq \frac{N(d,M)}{(1+|y|^2)^M}\sum_{|\alpha|\leq2M}\int_{\{\frac{1}{2}<|\xi|<2\}}\left|D_\xi^\alpha\left(|\xi|^\sigma\phi(\xi)\right)\right|d\xi
\\
&\leq N(d,M)(1+|y|^2)^{-M} \leq N(d,M)(1+|y|)^{-M},
\end{align*}
where $N$ is independent of $\sigma$.
Thus, for $M>d+2$ and $y\neq0$,
\begin{align} \label{eq11101320}
\sum_{j\geq1}|q_j(y)| &\leq N|y|^{-d-\sigma} \left(\sum_{\{j\geq1, 2^j|y|\leq1\}} (2^j|y|)^{d+\sigma} + \sum_{\{j\geq1, 2^j|y|>1\}}(2^j|y|)^{d+\sigma-M}\right) \nonumber
\\
&\leq N(d)|y|^{-d-\sigma}.
\end{align}

Let
\begin{equation*}
\psi(y):=|y|^{d+\sigma}q(y).
\end{equation*}
Then by \eqref{eq11101320} and $q\in \cS(\bR^d)$, one can see that $\psi\in\text{BUC}(\bR^d)$. 
Again by \eqref{eq11101320},
\begin{equation} \label{eq11101334}
\psi(2^jy)=|y|^{d+\sigma}q_j(y), \quad \sup_{y\in\bR^d}\sum_{j\geq1}|\psi(2^jy)|\leq N(d).
\end{equation}
Let $\delta\in(0,1)$ and take $a_*=a_*(d,\delta)>0$ so that $N(d)a_*\leq\delta$, where $N(d)$ is taken from \eqref{eq11101334}.
Define
\begin{equation} \label{eq11101338}
a:=a_*c_{d,\sigma}, \quad a_j:=\frac{a}{\sqrt{j}}, \quad j\geq1.
\end{equation}

Let $\beta_1:=3$. Choose $\gamma_m\geq\beta_m$ to be the smallest integer such that
\begin{equation*}
S_m:=\sum_{j=\beta_m}^{\gamma_m}a_j^2\geq m^2,
\end{equation*}
and let $\beta_{m+1}:=\gamma_m+3$.
The existence of $\gamma_m$ and $\beta_m$ is guaranteed due to $\sum_j a_j^2 = \infty$.
Moreover,
\begin{equation*}
\beta_m\geq\beta_1+3(m-1)=3m
\end{equation*}
and for $m\geq2$, we have
\begin{equation} \label{eq11101644}
\gamma_{m-1}<\beta_m-2=\gamma_{m-1}+1<\beta_m.
\end{equation}

Now we construct the desired kernel.
For $j\geq1$, let $b_j(x):=\cos(2^jx_1)$. Define
\begin{equation} \label{eq11101423}
K(x,y):=(2-\sigma)\left(1+\frac{1}{c_{d,\sigma}}\sum_{m\geq1}b_{\beta_m-2}(x)\sum_{j=\beta_m}^{\gamma_m}a_jb_j(x)\psi(2^jy)\right).
\end{equation}
By \eqref{eq11101334} and \eqref{eq11101338},
\begin{equation*}
\left|\frac{K(x,y)}{2-\sigma}-1\right|\leq\frac{a}{c_{d,\sigma}}\sum_{j\geq1}|\psi(2^jy)|\leq a_*N(d) \leq \delta,
\end{equation*}
which yields
\begin{equation*}
(1-\delta)(2-\sigma)\leq K(x,y)\leq(1+\delta)(2-\sigma).
\end{equation*}

For $n\in\bN$, define
\begin{equation*}
K_n(x,y):=(2-\sigma) \left(1+\frac{1}{c_{d,\sigma}} \sum_{m=1}^{n} b_{\beta_m-2}(x) \sum_{j=\beta_m}^{\gamma_m} a_jb_j(x) \psi(2^jy)\right).
\end{equation*}
Then by \eqref{eq11101334},
\begin{equation*}
\left\|\frac{K-K_n}{2-\sigma}\right\|_{L_\infty(\bR^{2d})} \leq N(d,\sigma) \sup_{j\geq\beta_{n+1}}a_j \to0
\end{equation*}
as $n\to\infty$.
Since $K_n\in \text{BUC}(\bR^d\times\bR^d)$, one can deduce that $K\in\text{BUC}(\bR^d\times\bR^d)$.
Since $\psi$ is radial, $K(x,\cdot)$ is radial for every $x\in\bR^d$.
In particular, the cancellation condition \eqref{eq_cancel} holds when $\sigma=1$. Therefore, $K$ satisfies Assumptions \ref{as_bound} and \ref{as_nonde} with $\Lambda=1+\delta$ and $\nu=(1-\delta)|\bS^{d-1}|/d$.

For a bounded smooth function $v$, define $T_jv:=q_j*v$.
Since $q_j$ is even and has zero integral,
\begin{equation*}
T_jv(x) = \int_{\bR^d} \delta_yv(x)q_j(y) dy.
\end{equation*}
In particular,
\begin{equation} \label{eq11101440}
T_jb_k = \hat{q}_j(2^ke_1) b_k = 2^{k\sigma} \phi(2^{k-j}e_1)b_k = 1_{\{j=k\}} 2^{k\sigma} b_k.
\end{equation}
Here, the last equality follows from \eqref{eq11101439}.

For a bounded smooth function $v$ with bounded derivatives, due to \eqref{intro_op} and \eqref{eq11101423},
\begin{equation} \label{eq11101441}
L_Kv = \frac{2-\sigma}{c_{d,\sigma}} \left(-(-\Delta)^{\sigma/2}v + \sum_{m\geq1} b_{\beta_m-2} \sum_{j=\beta_m}^{\gamma_m} a_jb_jT_jv\right).
\end{equation}
Here, the interchange of the sum and the integral in \eqref{intro_op} is justified by \eqref{eq11101320} and
\begin{equation*}
|v(x+y)+v(x-y)-2v(x)| \leq N(v) \min\{|y|^2,1\},
\end{equation*}
together with the symmetry of the kernels.

\subsection{Periodic test functions}

\begin{lemma}
Let $p\in(1,\infty)$ and let $\beta,\gamma$ be integers with $0\leq\beta\leq\gamma$. Then for any real numbers $t_\beta,\ldots,t_\gamma$,
\begin{equation} \label{eq11101525}
\left(\int_0^{2\pi} \left|\sum_{j=\beta}^{\gamma}t_j\cos(2^jt)\right|^p dt\right)^{1/p} \leq N(p) \left(\sum_{j=\beta}^{\gamma} |t_j|^2\right)^{1/2}.
\end{equation}
\end{lemma}

\begin{proof}
Let $F(t):=\sum_{j=\beta}^{\gamma}t_j\cos(2^jt)$ and $C^2:=\sum_{j=\beta}^{\gamma}t_j^2$ with $C>0$.

Note that
\begin{equation*}
\cos\left(2^\beta(t+\pi2^{-\beta})\right) = \cos(2^\beta t+\pi) = -\cos(2^\beta t)
\end{equation*}
and for $j=\beta+1,\ldots,\gamma$,
\begin{align*}
\cos\left(2^j(t+\pi2^{-\beta})\right) = \cos\left(2^jt+\pi2^{j-\beta}\right) =\cos(2^jt).
\end{align*}
This means that
\begin{equation} \label{eq11101538}
F(t+\pi2^{-\beta}) = -t_\beta \cos(2^\beta t) + \sum_{j=\beta+1}^{\gamma} t_j \cos(2^jt).
\end{equation}
Since
\begin{equation} \label{eq11101552}
\int_0^{2\pi}|F(t)|^p dt \leq N(p) C^p \int_0^{2\pi} \cosh\left(\frac{F(t)}{C}\right) dt,
\end{equation}
we aim to estimate the right-hand side of this inequality.
Using \eqref{eq11101538} and the identity $\cosh(a+b)+\cosh(a-b)=2\cosh(a)\cosh(b)$,
\begin{align*}
\int_0^{2\pi} \cosh\left(\frac{F(t)}{C}\right) dt &= \frac{1}{2} \int_0^{2\pi} \left(\cosh\left(\frac{F(t)}{C}\right) + \cosh\left(\frac{F(t+\pi2^{-\beta})}{C}\right)\right) dt
\\
&= \int_0^{2\pi} \cosh\left(\frac{1}{C} \sum_{j=\beta+1}^{\gamma} t_j\cos(2^jt)\right) \cosh\left(\frac{1}{C}t_\beta\cos(2^\beta t)\right) dt
\\
&\leq e^{\frac{t_\beta^2}{2C^2}} \int_0^{2\pi} \cosh\left(\frac{1}{C} \sum_{j=\beta+1}^{\gamma} t_j\cos(2^jt)\right) dt,
\end{align*}
where the last inequality follows from $\cosh z\leq e^{z^2/2}$ and $|\cos(2^\beta t)|\leq1$.
Repeating this argument,
\begin{equation*}
\frac{1}{2\pi} \int_0^{2\pi} \cosh\left(\frac{F(t)}{C}\right) dt \leq \exp\left(\frac{1}{2C^2} \sum_{j=\beta}^{\gamma}t_j^2\right) = e^{1/2}.
\end{equation*}
This and \eqref{eq11101552} yield \eqref{eq11101525}.
The lemma is proved.
\end{proof}

\begin{corollary} \label{cor11101719}
Let $d\geq1$, $p\in(1,\infty)$, and let $\beta,\gamma$ be integers with $0\leq\beta\leq\gamma$. Let $b_j(x)=\cos(2^jx_1)$. Then for any real numbers $t_\beta,\ldots,t_\gamma$,
\begin{equation*}
\left\|\sum_{j=\beta}^{\gamma}t_jb_j\right\|_{L_p(B_2)} \leq N(d,p) \left(\sum_{j=\beta}^{\gamma}|t_j|^2\right)^{1/2}.
\end{equation*}
\end{corollary}

\begin{proof}
Using $B_2\subset[-2,2]^d$ and \eqref{eq11101525}, we obtain
\begin{align*}
\left\|\sum_{j=\beta}^{\gamma}t_jb_j\right\|_{L_p(B_2)}^p &\leq 4^{d-1}\int_{-2}^{2} \left|\sum_{j=\beta}^{\gamma}t_j\cos(2^jt)\right|^p dt
\\
&\leq 4^{d-1} \int_0^{2\pi} \left|\sum_{j=\beta}^{\gamma}t_j\cos(2^jt)\right|^p dt \leq N(d,p) \left(\sum_{j=\beta}^{\gamma} |t_j|^2\right)^{p/2}.
\end{align*}
The corollary is proved.
\end{proof}

For each $m\geq1$, let $\rho_m := 2^{\beta_m-2}$ and define
\begin{align} \label{eq11101612}
h_m(x) &:= \frac{2}{S_m} \sum_{j=\beta_m}^{\gamma_m} a_j2^{-j\sigma}b_j(x), \nonumber
\\
u_m(x) &:= \rho_m^{-\sigma} b_{\beta_m-2}(x) + h_m(x).
\end{align}
Then they are smooth and have bounded derivatives of every order.

\begin{lemma}
Let $h_m$ and $u_m$ be defined by \eqref{eq11101612}. Then we have
\begin{align} \label{eq11101637}
(-\Delta)^{\sigma/2}h_m &= \frac{2}{S_m} \sum_{j=\beta_m}^{\gamma_m} a_jb_j, \nonumber
\\
(-\Delta)^{\sigma/2}u_m &= b_{\beta_m-2}+(-\Delta)^{\sigma/2}h_m.
\end{align}
Moreover,
\begin{equation} \label{eq11101639}
L_K(\rho_m^{-\sigma}b_{\beta_m-2}) = -\frac{2-\sigma}{c_{d,\sigma}} b_{\beta_m-2},
\end{equation}
and
\begin{align} \label{eq11101641}
\frac{c_{d,\sigma}}{2-\sigma} L_Ku_m = -\frac{2}{S_m} \sum_{j=\beta_m}^{\gamma_m} a_jb_j + \frac{b_{\beta_m-2}}{S_m} \sum_{j=\beta_m}^{\gamma_m} a_j^2b_{j+1}.
\end{align}
\end{lemma}

\begin{proof}
We first compute the fractional Laplacian of $b_j$ directly.
Recall that
\begin{equation} \label{eq1110400}
-(-\Delta)^{\sigma/2}= L_{c_{d,\sigma}},
\end{equation}
where $c_{d,\sigma}:=\frac{\sigma2^{\sigma-1}\Gamma((d+\sigma)/2)}{\pi^{d/2}\Gamma(1-\sigma/2)}$.
By using \eqref{eq1110400} and the identity
\begin{equation*}
2b_j(x)-b_j(x+y)-b_j(x-y)=2b_j(x)(1-\cos(2^jy_1)),
\end{equation*}
we have
\begin{align*}
(-\Delta)^{\sigma/2} b_j(x) &= \frac{c_{d,\sigma}}{2} \int_{\bR^d}\frac{2b_j(x)-b_j(x+y)-b_j(x-y)}{|y|^{d+\sigma}} dy
\\
&= c_{d,\sigma}b_j(x) \int_{\bR^d}\frac{1-\cos(2^jy_1)}{|y|^{d+\sigma}} dy
\\
&= 2^{j\sigma}c_{d,\sigma}b_j(x) \int_{\bR^d}\frac{1-\cos{z_1}}{|z|^{d+\sigma}} dz =2^{j\sigma}b_j(x).
\end{align*}
Here, for the last equality, we used
\begin{equation*}
c_{d,\sigma} \int_{\bR^d}\frac{1-\cos{z_1}}{|z|^{d+\sigma}}dz = 1.
\end{equation*}
Consequently,
\begin{align*}
(-\Delta)^{\sigma/2}h_m = \frac{2}{S_m} \sum_{j=\beta_m}^{\gamma_m} a_j2^{-j\sigma}(-\Delta)^{\sigma/2}b_j = \frac{2}{S_m} \sum_{j=\beta_m}^{\gamma_m} a_jb_j,
\end{align*}
and
\begin{align*}
(-\Delta)^{\sigma/2}u_m = \rho_m^{-\sigma} 2^{(\beta_m-2)\sigma} b_{\beta_m-2} + (-\Delta)^{\sigma/2} h_m = b_{\beta_m-2} + (-\Delta)^{\sigma/2} h_m,
\end{align*}
which proves \eqref{eq11101637}.

By \eqref{eq11101644} and \eqref{eq11101440}, if $\beta_n\leq j\leq\gamma_n$ for some $n\in\bN$, then
\begin{equation*}
T_j(\rho_m^{-\sigma}b_{\beta_m-2})=0.
\end{equation*}
This and \eqref{eq11101441} lead to \eqref{eq11101639}.

For $j\geq1$, due to \eqref{eq11101440}, we also have
\begin{equation*}
T_ju_m = 1_{\{j=\beta_m-2\}}b_{\beta_m-2} + \frac{2}{S_m}1_{\{\beta_m\leq j\leq\gamma_m\}} a_jb_j.
\end{equation*}
 Using \eqref{eq11101644} and $2b_j^2=1+b_{j+1}$, we obtain
\begin{align*}
\sum_{n\geq1} b_{\beta_n-2} \sum_{j=\beta_n}^{\gamma_n} a_jb_jT_ju_m &= \frac{2b_{\beta_m-2}}{S_m} \sum_{j=\beta_m}^{\gamma_m} a_j^2b_j^2
\\
&= \frac{b_{\beta_m-2}}{S_m} \sum_{j=\beta_m}^{\gamma_m} a_j^2(1+b_{j+1}) = b_{\beta_m-2} + \frac{b_{\beta_m-2}}{S_m} \sum_{j=\beta_m}^{\gamma_m} a_j^2b_{j+1}.
\end{align*}
This together with \eqref{eq11101441} and \eqref{eq11101637} yields \eqref{eq11101641}.
The lemma is proved.
\end{proof}

For $\alpha\in(0,1)$, we write
\begin{equation*}
\|v\|_{C^\alpha} := \|v\|_{L_\infty}+[v]_{C^\alpha}, \quad [v]_{C^\alpha} := \sup_{x\neq y} \frac{|v(x)-v(y)|}{|x-y|^\alpha}.
\end{equation*}

\begin{lemma}
Let $p\in(1,\infty)$, $\sigma\in(0,2)$, and $h_m$ and $u_m$ be defined by \eqref{eq11101612}. Then we have
\begin{equation} \label{eq11101718}
\|L_Ku_m\|_{L_p(B_2)}\leq NS_m^{-1/2},
\end{equation}
\begin{equation} \label{eq11101719}
\|(-\Delta)^{\sigma/2}h_m\|_{L_p(B_2)}\leq NS_m^{-1/2},
\end{equation}
and
\begin{equation} \label{eq11101720}
\|u_m\|_{C^{\sigma/2}}+\|h_m\|_{C^{\sigma/2}}\leq N2^{-\beta_m\sigma/2}(1+S_m^{-1/2}),
\end{equation}
where $N=N(d,p,\widetilde{\sigma},\delta)$ and $\widetilde{\sigma}$ is given by \eqref{eq_sigma}.
Moreover,
\begin{equation} \label{eq11101721}
\liminf_{m\to\infty}\|b_{\beta_m-2}\|_{L_p(B_1)}\geq c(d,p)>0.
\end{equation}
\end{lemma}

\begin{proof}
We first consider \eqref{eq11101719}.
By \eqref{eq11101637} and Corollary \ref{cor11101719} with $\beta=\beta_m$, $\gamma=\gamma_m$, and $t_j=a_j$ for $j=\beta_m,\ldots,\gamma_m$,
\begin{align*}
\|(-\Delta)^{\sigma/2}h_m\|_{L_p(B_2)} &= \frac{2}{S_m}\left\|\sum_{j=\beta_m}^{\gamma_m}a_jb_j\right\|_{L_p(B_2)} \leq \frac{N(d,p)}{S_m}\left(\sum_{j=\beta_m}^{\gamma_m}a_j^2\right)^{1/2}
\\
&=\frac{N(d,p)}{S_m}S_m^{1/2} = N(d,p)S_m^{-1/2}.
\end{align*}
 
 Next, we prove \eqref{eq11101718}.
By \eqref{eq11101338}, \eqref{eq11101641}, and Corollary \ref{cor11101719},
\begin{align*}
\|L_Ku_m\|_{L_p(B_2)} &\leq \frac{2-\sigma}{c_{d,\sigma}S_m} \left(2\left\|\sum_{j=\beta_m}^{\gamma_m}a_jb_j\right\|_{L_p(B_2)} + \left\|\sum_{j=\beta_m}^{\gamma_m}a_j^2b_{j+1}\right\|_{L_p(B_2)}\right)
\\
&\leq \frac{N(d,p)(2-\sigma)}{c_{d,\sigma}S_m} \left(\left(\sum_{j=\beta_m}^{\gamma_m}a_j^2\right)^{1/2} + \left(\sum_{j=\beta_m}^{\gamma_m}a_j^4\right)^{1/2}\right)
\\
&\leq \frac{N}{S_m} \left(S_m^{1/2}+aS_m^{1/2}\right) = N(1+a)S_m^{-1/2} \leq N S_m^{-1/2}.
\end{align*}
Thus, we have \eqref{eq11101718}.

We next deal with \eqref{eq11101720}. Note that for $c\geq1$ and $\alpha\in(0,1)$,
\begin{equation*}
|\cos(c x_1)-\cos(c z_1)|\leq\min\{2,c|x-z|\}\leq 2c^\alpha|x-z|^\alpha.
\end{equation*}
This together with $|\cos(c x_1)|\leq1$ leads to
\begin{equation*}
\|\cos(c x_1)\|_{C^\alpha} \leq 1 + 2c^\alpha \leq 3 c^\alpha.
\end{equation*}
By applying the Cauchy-Schwarz inequality, we obtain
\begin{align*}
\|h_m\|_{C^{\sigma/2}} &\leq \frac{6}{S_m}\sum_{j=\beta_m}^{\gamma_m} a_j2^{-j\sigma/2} \leq 6S_m^{-1/2} \left(\sum_{j\geq\beta_m} 2^{-j\sigma}\right)^{1/2}
\\
&=\frac{6S_m^{-1/2} 2^{-\beta_m\sigma/2}}{\sqrt{1-2^{-\sigma}}} \leq N(\widetilde{\sigma}) S_m^{-1/2}2^{-\beta_m\sigma/2}.
\end{align*}
Since $\rho_m^{-\sigma/2} = 2^\sigma2^{-\beta_m\sigma/2} \leq 4\cdot2^{-\beta_m\sigma/2}$, it follows from \eqref{eq11101612} that
\begin{align*}
\|u_m\|_{C^{\sigma/2}}+\|h_m\|_{C^{\sigma/2}} \leq 3\rho_m^{-\sigma/2}+2\|h_m\|_{C^{\sigma/2}} \leq N(\widetilde{\sigma})2^{-\beta_m\sigma/2}(1+S_m^{-1/2}),
\end{align*}
which proves \eqref{eq11101720}. 

Lastly, we consider \eqref{eq11101721}.
Let $l_0=\frac{1}{2\sqrt{d}}$ so that $(-l_0,l_0)^d\subset B_1$.
Then
\begin{align*}
\int_{B_1} |b_{\beta_m-2}(x)|^p dx &\geq(2l_0)^{d-1}\int_{-l_0}^{l_0}|\cos(\rho_m t)|^p dt =\frac{(2l_0)^{d-1}}{\rho_m} \int_{-l_0\rho_m}^{l_0\rho_m} |\cos{t}|^p dt.
\end{align*}
Let $n:=\left\lfloor\frac{l_0\rho_m}{2\pi}\right\rfloor$. Then $(-2\pi n,2\pi n)\subset(-l_0\rho_m,l_0\rho_m)$, which yields
\begin{align*}
\int_{B_1} |b_{\beta_m-2}(x)|^p dx &\geq\frac{(2l_0)^{d-1}}{\rho_m}\int_{-2\pi n}^{2\pi n}|\cos{t}|^p dt
\\
&= \frac{(2l_0)^{d-1}}{\rho_m}2n \int_0^{2\pi}|\cos{t}|^p dt
\\
&\geq (2l_0)^{d-1} \left(2l_0-\frac{4\pi}{\rho_m}\right) \frac{1}{2\pi} \int_0^{2\pi} |\cos{t}|^p dt,
\end{align*}
where the last inequality follows from $n\geq\frac{l_0\rho_m}{2\pi}-1$.
Since $\rho_m \geq 2^{3m-2}$, there is $m_0=m_0(d)$ such that $4\pi/\rho_m \leq l_0$ for $m\geq m_0$. Thus,
\begin{equation*}
\|b_{\beta_m-2}\|_{L_p(B_1)}^p \geq \frac{(2l_0)^{d-1}l_0}{2\pi}\int_0^{2\pi} |\cos{t}|^p dt>0, \quad m\geq m_0,
\end{equation*}
which proves \eqref{eq11101721}.
The lemma is proved.
\end{proof}

\subsection{Proof of the counterexamples}

\begin{lemma}
Let $\sigma\in(0,2)$ and $p\in(1,\infty)$. Let $\zeta\in C_c^\infty(B_2)$ such that $0\leq\zeta\leq1$ and $\zeta=1$ on $B_1$.
Suppose that $K$ satisfies Assumption \ref{as_bound} and
\begin{equation*}
K(x,y)=K(x,-y).
\end{equation*}
Then for any bounded smooth function $v$ with bounded derivatives, we have
\begin{equation} \label{eq11102021}
\|[L_K,\zeta]v\|_{L_p} \leq N\Lambda \|v\|_{C^{\sigma/2}},
\end{equation}
where $N=N(d,p,\widetilde{\sigma},\zeta)$ and $\widetilde{\sigma}$ is given by \eqref{eq_sigma}.
Moreover,
\begin{equation} \label{eq11102022}
\|[(-\Delta)^{\sigma/2},\zeta]v\|_{L_p} \leq N\|v\|_{C^{\sigma/2}}.
\end{equation}
\end{lemma}

\begin{proof}
Due to \eqref{eq1110400}, we only need to prove \eqref{eq11102021}.
Note that
\begin{align} \label{eq11102029}
[L_K,\zeta]v(x) = v(x)L_K\zeta(x) + \int_{\bR^d} (\zeta(x+y)-\zeta(x))(v(x+y)-v(x)) K(x,y) \frac{dy}{|y|^{d+\sigma}}.
\end{align}
Due to the symmetry of $K$,
\begin{equation*}
L_K\zeta(x) = \frac{1}{2} \int_{\bR^d} (\zeta(x+y)+\zeta(x-y)-2\zeta(x))K(x,y) \frac{dy}{|y|^{d+\sigma}}.
\end{equation*}

We first estimate $[L_K,\zeta]v$ when $x\in B_4$. For the first term on the right-hand side of \eqref{eq11102029}, we have
\begin{align} \label{eq11102039}
\|L_K\zeta\|_{L_\infty}&\leq N(\zeta) \Lambda(2-\sigma) \left(\int_{B_1} |y|^{2-d-\sigma} dy + \int_{B_1^c} |y|^{-d-\sigma} dy\right) \nonumber
\\
&= N(\zeta) \Lambda (2-\sigma) \left(\frac{1}{2-\sigma} + \frac{1}{\sigma}\right) \leq N(\zeta) \Lambda.
\end{align}

For the second term of the right-hand side of \eqref{eq11102029},
\begin{align*}
&\int_{\bR^d} |\zeta(x+y)-\zeta(x)||v(x+y)-v(x)|K(x,y) \frac{dy}{|y|^{d+\sigma}}
\\
&\leq N(\zeta) \Lambda (2-\sigma) \left([v]_{C^{\sigma/2}}\int_{B_1} |y|^{1-d-\sigma/2} dy + \|v\|_{L_\infty} \int_{B_1^c} |y|^{-d-\sigma} dy\right)
\\
&= N(\zeta) \Lambda (2-\sigma) \left(\frac{[v]_{C^{\sigma/2}}}{1-\sigma/2} + \frac{\|v\|_{L_\infty}}{\sigma}\right) \leq N\Lambda \|v\|_{C^{\sigma/2}}.
\end{align*}
This together with \eqref{eq11102039} leads to
\begin{equation} \label{eq11102046}
\|[L_K,\zeta]v\|_{L_p(B_4)} \leq N\Lambda \|v\|_{C^{\sigma/2}}.
\end{equation}

Now we consider $[L_K,\zeta]v$ when $|x|\geq4$.
In this case,
\begin{equation*}
[L_K,\zeta]v(x) = \int_{B_2}\zeta(z)v(z)K(x,z-x) \frac{dz}{|z-x|^{d+\sigma}}.
\end{equation*}
Since $|z-x|\geq|x|/2$ for $z\in B_2$, we have
\begin{equation*}
|[L_K,\zeta]v(x)| \leq N \Lambda (2-\sigma) |x|^{-d-\sigma}\|v\|_{L_\infty}.
\end{equation*}
Hence,
\begin{align} \label{eq11102047}
\|[L_K,\zeta]v\|_{L_p(\bR^d\setminus B_4)}^p &\leq N\Lambda^p \|v\|_{L_\infty}^p \int_{\bR^d\setminus B_4} |y|^{-p(d+\sigma)} dy \nonumber
\\
&=N\Lambda^p \|v\|_{L_\infty}^p \frac{4^{d-p(d+\sigma)}}{p(d+\sigma)-d} \leq N\Lambda^p\|v\|_{L_\infty}^p.
\end{align}
Thus, \eqref{eq11102021} follows from \eqref{eq11102046} and \eqref{eq11102047}.
The lemma is proved.
\end{proof}

\begin{proof}[Proof of Theorem \ref{thm_counter}]

Let $K$ be the kernel defined in \eqref{eq11101423}.
By the construction in Section \ref{sec_cons}, $K\in\text{BUC}(\bR^d\times\bR^d)$ and satisfies Assumptions
\ref{as_bound} and \ref{as_nonde}.
We consider functions $\zeta h_m$ and $\zeta u_m$ in $C_c^\infty(B_2)$, where $\zeta\in C_c^\infty(B_2)$ satisfies $0\leq\zeta\leq1$ and $\zeta=1$ on $B_1$.

$(i)$
We prove that \eqref{eq_buc_bound} fails by showing that
\begin{equation} \label{eq11102114}
\lim_{m\to\infty} \|\zeta h_m\|_{H_p^\sigma}=0 \text{ and } \liminf_{m\to\infty} \|L_K(\zeta h_m)\|_{L_p} > 0.
\end{equation}
Since
\begin{equation*}
(-\Delta)^{\sigma/2}(\zeta h_m)=\zeta(-\Delta)^{\sigma/2}h_m+[(-\Delta)^{\sigma/2},\zeta]h_m,
\end{equation*}
by \eqref{eq10310040}, \eqref{eq11101719}, \eqref{eq11101720}, and \eqref{eq11102022}, we obtain
\begin{align*}
\|\zeta h_m\|_{H_p^\sigma} &\leq N\left(\|\zeta h_m\|_{L_p} + \|(-\Delta)^{\sigma/2}(\zeta h_m)\|_{L_p}\right)
\\
&\leq N \left(\|h_m\|_{L_\infty}+\|(-\Delta)^{\sigma/2}h_m\|_{L_p(B_2)}+\|[(-\Delta)^{\sigma/2},\zeta]h_m\|_{L_p}\right)
\\
&\leq N\left(\|h_m\|_{C^{\sigma/2}}+\|(-\Delta)^{\sigma/2}h_m\|_{L_p(B_2)}\right)\to0
\end{align*}
as $m\to\infty$.

By \eqref{eq11101612} and \eqref{eq11101639},
\begin{equation*}
L_Kh_m = L_Ku_m+\frac{2-\sigma}{c_{d,\sigma}}b_{\beta_m-2}.
\end{equation*}
Therefore,
\begin{equation*}
L_K(\zeta h_m)-\frac{2-\sigma}{c_{d,\sigma}}\zeta b_{\beta_m-2} = \zeta L_Ku_m+[L_K,\zeta]h_m.
\end{equation*}
By \eqref{eq11102021} with $\Lambda=1+\delta$, \eqref{eq11101718}, and \eqref{eq11101720}, we obtain
\begin{equation*}
\left\|L_K(\zeta h_m)-\frac{2-\sigma}{c_{d,\sigma}}\zeta b_{\beta_m-2}\right\|_{L_p} \leq \|L_Ku_m\|_{L_p(B_2)} + N \|h_m\|_{C^{\sigma/2}}\to0
\end{equation*}
as $m\to\infty$.
Moreover, due to \eqref{eq11101721},
\begin{equation*}
\liminf_{m\to\infty}\|L_K(\zeta h_m)\|_{L_p} \geq \frac{2-\sigma}{c_{d,\sigma}} \liminf_{m\to\infty} \|b_{\beta_m-2}\|_{L_p(B_1)} > 0,
\end{equation*}
which proves \eqref{eq11102114}. Applying \eqref{eq_buc_bound} to $\zeta h_m$ would give a contradiction.

$(ii)$
We first prove that
\begin{align} \label{eq11102115}
&\lim_{m\to\infty} \left(\|\zeta u_m\|_{L_p}+\|L_K(\zeta u_m)\|_{L_p}\right) =0, \nonumber
\\
&\liminf_{m\to\infty}\|(-\Delta)^{\sigma/2}(\zeta u_m)\|_{L_p} > 0.
\end{align}
Since
\begin{equation*}
L_K(\zeta u_m) = \zeta L_Ku_m+[L_K,\zeta]u_m,
\end{equation*}
by \eqref{eq11101720}, \eqref{eq11102021}, and \eqref{eq11101718},
\begin{equation*}
\|\zeta u_m\|_{L_p}+\|L_K(\zeta u_m)\|_{L_p}\leq N\|u_m\|_{C^{\sigma/2}}+\|L_Ku_m\|_{L_p(B_2)}\to0
\end{equation*}
as $m\to\infty$.

On the other hand, due to \eqref{eq11101637}, we have
\begin{equation*}
(-\Delta)^{\sigma/2}(\zeta u_m)-\zeta b_{\beta_m-2} = \zeta(-\Delta)^{\sigma/2}h_m+[(-\Delta)^{\sigma/2},\zeta]u_m.
\end{equation*}
Thus, by \eqref{eq11101719}, \eqref{eq11101720}, and \eqref{eq11102022},
\begin{equation*}
\|(-\Delta)^{\sigma/2}(\zeta u_m)-\zeta b_{\beta_m-2}\|_{L_p}\leq\|(-\Delta)^{\sigma/2}h_m\|_{L_p(B_2)}+N\|u_m\|_{C^{\sigma/2}}\to0
\end{equation*}
as $m\to\infty$.
Therefore, this and \eqref{eq11101721} imply
\begin{equation*}
\liminf_{m\to\infty} \|(-\Delta)^{\sigma/2}(\zeta u_m)\|_{L_p} \geq \liminf_{m\to\infty} \|b_{\beta_m-2}\|_{L_p(B_1)}>0.
\end{equation*}
Thus, \eqref{eq11102115} is proved.

Now we conclude that \eqref{eq_buc_est} fails.
For any $\lambda\geq0$,
\begin{align*}
\|L_K(\zeta u_m)-\lambda\zeta u_m\|_{L_p} + \|\zeta u_m\|_{L_p} \leq\|L_K(\zeta u_m)\|_{L_p}+(1+\lambda)\|\zeta u_m\|_{L_p} \to 0
\end{align*}
as $m\to\infty$.
This and \eqref{eq11102115} yield a contradiction to \eqref{eq_buc_est}.
The theorem is proved.
\end{proof}


\section*{Acknowledgement} \pdfbookmark{Acknowledgement}{Acknowledgement}

Section \ref{sec_counter} was reorganized from \cite{L26} with the assistance of ChatGPT 6 Astra and revised by the author. The author takes full responsibility for the content of this section.




\end{document}